\documentclass[12pt]{article}
\usepackage{subcaption}
\usepackage[a4paper, margin=2.5cm]{geometry}
\usepackage{setspace}
\usepackage{amsmath}               
\usepackage{amsthm}
\usepackage{amssymb}
\usepackage[title]{appendix}
\usepackage{hyperref}
\usepackage{cleveref}
\usepackage{xcolor}
\usepackage{graphicx}
\usepackage{rotating}
\newcommand*\rot{\rotatebox{90}}
\usepackage{comment}
\usepackage{bbm}
\usepackage{microtype}
\usepackage{authblk}
\usepackage{calc}
\theoremstyle{plain}
\newtheorem{theorem}{Theorem}[section]
\newtheorem*{theorem*}{Theorem}
\newtheorem{lemma}[theorem]{Lemma}
\newtheorem{claim}[theorem]{Claim}

\newtheorem*{assumption*}{Assumptions}
\newtheorem*{assumption}{Assumptions}

\newtheorem{definition}[theorem]{Definition}

\newcommand{\Var}{\ensuremath{\textrm{Var}}}

\newcommand{\norm}[1]{\left \lVert #1 \right \rVert}

\newcommand{\argmax}{\operatornamewithlimits{argmax}}

\usepackage{color}
\usepackage{xparse}

\providecommand{\keywords}[1]
{
  \small	
  \textbf{\textit{Keywords---}} #1
}

\makeatletter
\let\@fnsymbol\@arabic
\makeatother
\date{}

\title{Object detection under the linear subspace  model with application to cryo-EM images}

\author{Amitay Eldar $^{*,}$\thanks{Department of Applied Mathematics, Tel-Aviv University, Israel. \textsf{amitayeldar@tauex.tau.ac.il}}
	\qquad Keren Mor Waknin $^{*,}$\thanks{Department of Applied Mathematics, Tel-Aviv University, Israel. \textsf{kerenmor@tauex.tau.ac.il}} \qquad
	Samuel Davenport \thanks{Division of Biostatistics, University of California, San Diego, USA. \textsf{sdavenport@health.ucsd.edu}}	
	\qquad Tamir Bendory \thanks{School of Electrical Engineering, Tel Aviv University, Israel. \textsf{bendory@tauex.tau.ac.il}}
	\qquad Armin Schwartzman\thanks{\parbox{\textwidth-1cm}{Hal{\i}c{\i}o\u{g}lu Data Science Institute and Division of Biostatistics, University of California, San Diego, USA.  \textsf{armins@ucsd.edu}}\smallskip}
	\qquad Yoel Shkolnisky\thanks{Department of Applied Mathematics, Tel-Aviv University, Israel. \textsf{yoelsh@tauex.tau.ac.il}}}

\begin{document}
	\maketitle

\begin{abstract}
    Detecting multiple unknown objects in noisy data is a key problem in many scientific fields, such as electron microscopy imaging. A common model for the unknown objects is the linear subspace model, which assumes that the objects can be expanded in some known basis (such as the Fourier basis). In this paper, we develop an object detection algorithm that under the linear subspace model is asymptotically guaranteed to detect all objects, while controlling the family wise error rate or the false discovery rate. Numerical simulations show that the algorithm  also controls the error rate with high power in the non-asymptotic regime, even in highly challenging regimes.  We apply the proposed algorithm to experimental electron microscopy data set, and show that it outperforms existing standard software.
\end{abstract} \hspace{10pt}
\def\thefootnote{*}\footnotetext{These authors contributed equally}
\section{Introduction}\label{sec:introduction}

Detecting objects buried in noise is a fundamental problem in signal and image analysis, with applications in a wide variety of  scientific areas, such as neuroimaging~\cite{Genovese2002,Worsley1996,Taylor2007,davenport2022confidence}, electron microscopy~\cite{eldar2020klt,Bepler2019}, fluorescence microscopy~\cite{Geisler2007}, and astronomy~\cite{Brutti2005}. Nonetheless, in most cases there is no known optimal procedure for detecting the objects, nor algorithms with performance guarantees, and various problem-dependent heuristics are being used.

Many object detection algorithms were developed over the years. A comprehensive survey of object detection methods (about 50) is given in~\cite{ehret2019image}, focusing on anomaly detection (the problem which we call ``object detection'' is referred to in the literature by various names). In this problem, the image contains background (noise) and anomalies (objects) spread at unknown locations. The survey~\cite{ehret2019image} categorizes the different methods to probabilistic methods, distance-based methods, reconstruction-based methods, domain-based methods, and information-theoretic methods. Since we aim at algorithms with performance bounds, from our perspective, the methods in~\cite{ehret2019image} can be divided into two families. The first family includes algorithms with no theoretical guarantees whose results are tailored for specific data sets, see for example \cite{xie2007texems,perng2010novel,an2015variational,zontak2009defect,tsai1999automated}. The second family includes methods that are accompanied by theoretical guarantees that their family-wise error rate  (FWER~\cite{rupert2012simultaneous}) is controlled when the parameters of the algorithm are properly tuned~\cite{ehret2019image,davy2018reducing,aiger2010phase,zontak2010defect,grosjean2009contrario,boracchi2014novelty}. Importantly, none of these methods control the power and FWER simultaneously. Moreover, none of the algorithms provably controls the false discovery rate (FDR) \cite{benjamini1995controlling}.

Most recently,  \cite{marandon2022machine}  proposed an object detection algorithm that controls the FDR and maximizes power, assuming, among other things, that the observables share the same marginal distribution under the null hypothesis and exhibit exchangeability. However, as we aim to detect points close to an object's center, the null hypothesis in our model contains regions with objects, thus the observable does not share the same marginal distribution under the null hypotheses. Moreover, the test sample (referred to as the candidate set in our case~\ref{sec:algorithm}) does not satisfy exchangeability, preventing us from using the algorithm proposed by \cite{marandon2022machine}.

A  step towards deriving detection algorithms with performance guarantees for both the power and the error rate  (FDR or FWER) was carried out in~\cite{Adler2011}, under the rather restrictive model that the objects are one-dimensional, unimodal, and positive (namely, each object is positive and has only one maximum within its support). Under this model,~\cite{Adler2011} suggests a detection algorithm that asymptotically controls the FWER or the FDR, while its power tends asymptotically to one. The algorithm consists of smoothing the observed data followed by applying multiple hypothesis testing to the local maxima of the smoothed data. Specifically, the algorithm computes the $p$-values of the heights of the local maxima of the smoothed data, and then applies the Bonferroni (Bon)~\cite{rupert2012simultaneous} or Benjamini-Hochberg (BH)~\cite{benjamini1995controlling} procedures for thresholding the $p$-values. This provides asymptotic control of the FWER or the FDR with power tends to one, as the domain size and the norm of the objects grow. The algorithm relies on a closed-form expression for the distribution of the height of local maxima, which is unavailable in most practical cases. The approach of~\cite{Adler2011} has been extended in~\cite{Schwartzman2017} to higher dimensional Euclidean domains, yet the latter suffers from the same shortcomings as~\cite{Adler2011}. In particular, the objects to detect have to be unimodal. Moreover, both~\cite{Adler2011,Schwartzman2017}  only show that the peaks detected by their algorithms fall within the support of the objects, and are not necessarily close to the objects' centers (more on that issue below). 

In this paper, we design an object detection algorithm that overcomes the shortcomings of~\cite{Adler2011,Schwartzman2017}. First, the objects do not have to be unimodal but rather are spanned by an arbitrary orthogonal basis. This is what we call the linear subspace model. Second, we prove that with high probability, the estimated centers by the algorithm are close to the objects' centers, rather than just within the support of the objects. Third, as our setting requires a test statistic that is  Chi-Squared and not Gaussian, we derive an appropriate test function which allows to asymptotically control the error rate while achieving maximum power. The precise setting for our model is described in Section~\ref{sec:setup}. Our detection algorithm is described in Section~\ref{sec:algorithm}. The performance guarantees of the algorithm are provided in Section~\ref{sec:theorem}.

\begin{figure}
		\centering
			\centering
			\includegraphics[height=0.43\textwidth]{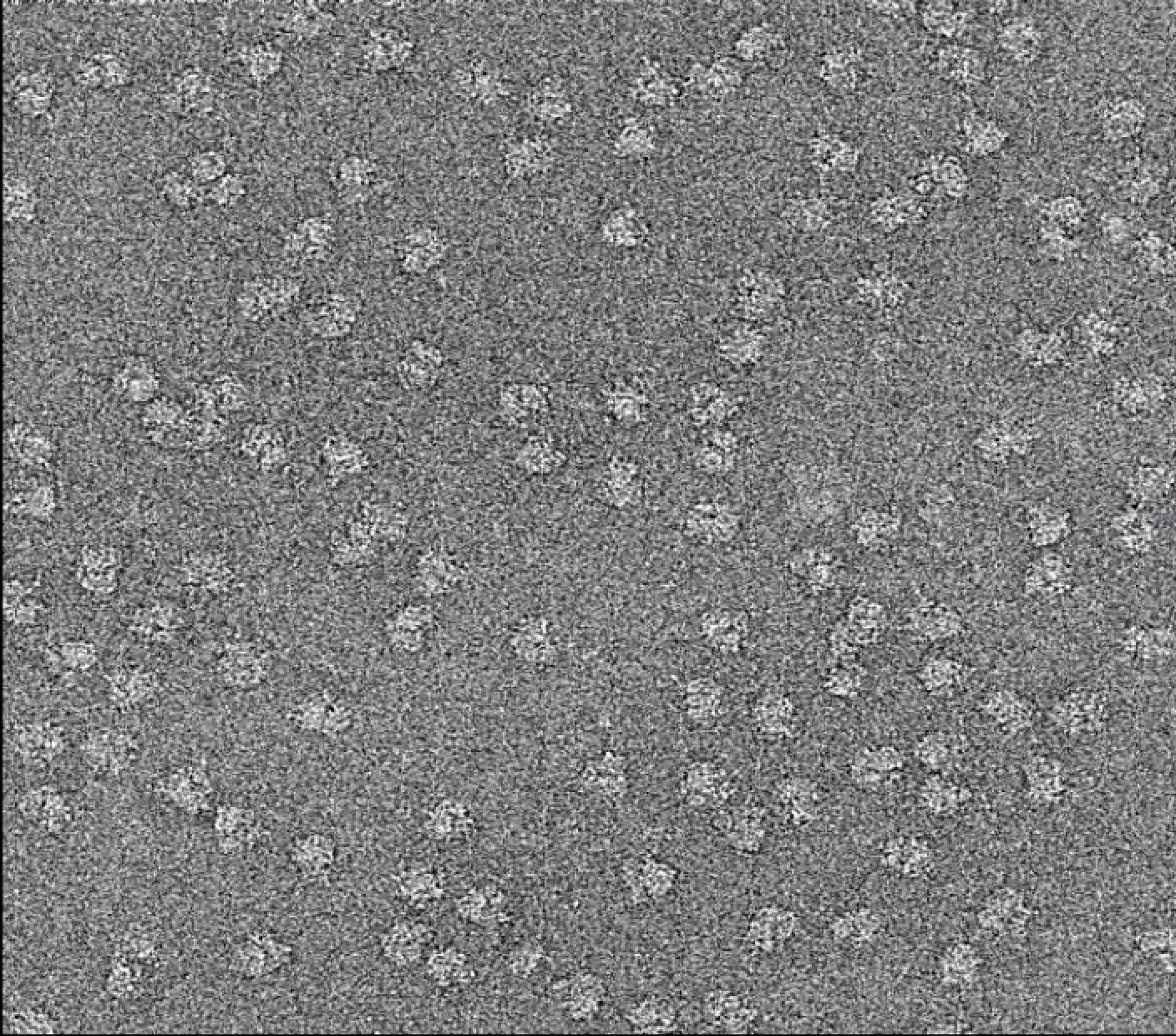}
			\includegraphics[height=0.43\textwidth]{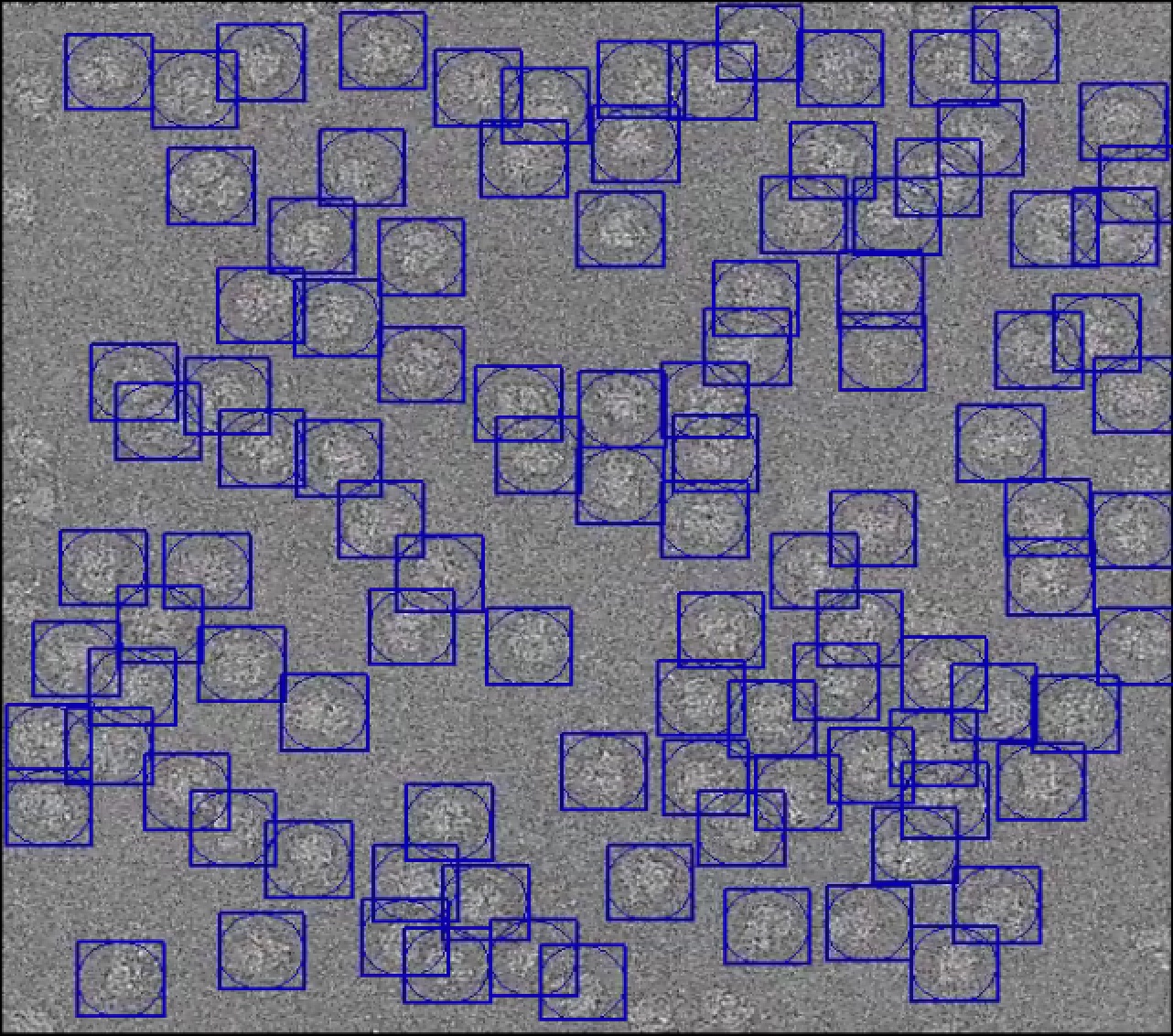}
		\caption{Example of a typical cryo-EM image (left panel) taken from \cite{10.7554/eLife.03080} and a manual detection output (right panel). }
            \label{fig: Cryo-EM}
	\end{figure}

A particular application of great scientific importance that motivates this work is the particle detection task in single-particle cryo-electron microscopy (cryo-EM), which is currently the leading method for resolving the three-dimensional structure of molecules~\cite{Lyumkis2019,NatureMethodEditorial}. In this method, a sample containing multiple copies of a molecule is rapidly frozen in a thin layer of ice, whereby the locations and orientations of the frozen molecules in the ice layer are unknown. An electron beam is transmitted through the frozen sample, resulting in a large, highly noisy, two-dimensional image containing two-dimensional smaller images of the molecules (technically, their tomographic projections) at unknown locations within the larger image. These two-dimensional small images are  typically referred to as 
particles. An example cryo-EM image together with manually annotated particles is shown in Figure~\ref{fig: Cryo-EM}. One of the first steps towards resolving the three-dimensional structure of the molecule is locating the particles in the image, a step called ``particle picking.'' Currently, all particle picking algorithms, e.g., \cite{relion3,EMAN,xmipp,Eldar2022} suffer from a high percentage of false positives (detecting ``non-particles'' as particles), thus requiring a time-consuming process of cleaning the data. The algorithm derived in this paper will allow to collect electron microscopy data sets with a low percentage of non-particle images, and will significantly speed up the processing of cryo-electron microscopy data. Moreover, it will eliminate the necessary manual intervention in current processing pipelines, which may introduce subjective bias into the analysis \cite{henderson2013avoiding,bendory2020single}.
  We will demonstrate in Section~\ref{sec: numerical results} the applicability of our algorithm to simulated and experimental cryo-EM data sets, and present improvement compared to one of the existing standard software. Based on our theoretical derivations (Section~\ref{sec:theorem}) and empirical insights (Section~\ref{sec: numerical results}), we discuss in Section~\ref{sec:future} the limitations and possible extensions of our algorithm.

To summarize, the contribution of our paper is five-fold. First, we introduce a general and flexible observation model for object detection tasks. Second, we devise an object detection algorithm accompanied by performance guarantees for its power and error rate.   In particular, we suggest  a  test statistic which is applicable to non-Gaussian data. Third, we introduce the notion of localization, which quantifies how accurately can we detect the centers of the objects. Fourth, we develop new proof techniques, which are essential in order to prove that the power of our algorithm approaches~$1$ while controlling the error rate. Finally, we apply our algorithm to experimental cryo-EM data sets, demonstrating improvement compared to existing algorithms.

The data analysis and all simulations were implemented in
Matlab. The code is available at  \url{https://github.com/ShkolniskyLab/object_detection_LSM/}.

\section{Theory}
\subsection{Problem setup}\label{sec:setup}
In this section, we introduce the observation model and define the object detection problem. 
We assume that the observable $y(t)$ is given by 
\begin{equation}\label{eq: the model}
	y(t)=x(t)+z(t), 
\end{equation}
where $x(t)$ consists of the objects to be  detected and $z(t)$ is a centered Gaussian stationary ergodic process on $\mathbb{R}^d$, representing the noise in the data. We assume that $y$ is observed on $\overline{C(0,L)}$, where 
\begin{equation}\label{eq:closed cube}
	\overline{C(t_0,L)} = \left \{t\in \mathbb{R}^d : \norm{t-t_0}_{\infty}\leq\frac{L}{2} \right \}
\end{equation}
is the closed hyper cube centered at $t_0$ with side length $L$. If the hyper cube in~\eqref{eq:closed cube} is open, we denote it by $C(t_{0},L)$. Additionally, suppose that $x(t)$ is given by 
\begin{equation}\label{eq:x}
x(t)=\sum_{i=1}^{N}x_i\left(t\right),\quad  x_i(t)= \sum_{j=1}^M a_{ij}\psi_j(t-\tau_i),
\end{equation}
where $N=N(L)$ is the (unknown) number of objects, $\left\{\psi_j\right\}_{j=1}^M$, are known functions, compactly supported on $C(0,B)$, $ C^\infty $ and orthonormal in $C(0,B)$, while the coefficients $ \left\{a_{ij}\right\}_{i=1,j=1}^{N,M} $ are unknown and depend on $L$. Note that each $ x_i $ in~\eqref{eq:x} is a deterministic function which is compactly supported within an open hypercube $C(\tau_i,B)\subset\mathbb{R}^d$ of a known side length~$B$, centered at an unknown location $\tau_i$. We term this model the \emph{linear subspace model}. Our goal is to estimate the objects' centers $\{\tau_i\}_{i=1}^N$ given $y(t)$, $B$, and $\{\psi_j\}_{j=1}^M$. 

Note that the proposed model is significantly more general than the model of~\cite{Adler2011,Schwartzman2017}, as the objects are not required to be unimodal nor positive, and we only assume that the objects are spanned by some known basis functions, such as the Fourier basis, orthogonal polynomials, PCA basis, etc. Moreover, in contrast to~\cite{Adler2011,Schwartzman2017}, we aim to estimate the unknown centers $\tau_{i}$, and not just detect some arbitrary point within the support of each object.

\subsection{Algorithm's outline}\label{sec:algorithm}

A typical object detection algorithm consists of two main steps: a)~applying a~\emph{score map (or a heat map)} on the observed data, which is a function that assigns a number to each position in the data;  b)~declaring all peaks in the score map larger than some \emph{threshold} as the centers of the detected objects. We consider each value of the score map as a proxy for the presence of an object at that location. Two fundamental questions are how to design the score map and how to set the threshold so that the detection algorithm has high power and controlled FWER or FDR.

We first introduce the algorithm step by step, and then discuss its different components.  
	\begin{enumerate}
		\item Compute the score map 
        \begin{align}\label{eq: smoothed process}
        S^y(t) =\sum_{j=1}^{M}(y \ast\tilde{\psi}_j )^2(t),
        \end{align}
        where $\tilde{\psi}_j(t)=\psi_j(-t)$ and $ \ast $ is a linear convolution over $\overline{C(0,L)}$.    
		\item Select candidate peaks. We choose the highest peak in $S^{y}$ to be our first candidate for an object's location. We then erase a box of size $r$ around this peak (the value of the parameter $r$ will be defined below). Then, we select the highest peak in the remaining data, and repeat this process until the remaining data is empty. Formally, we define a sequence of points
		\begin{equation*}
			t_1=\argmax_{\overline{C(0,L)}}S^y(t), \quad t_2=\argmax_{\overline{C(0,L)}\backslash C(t_1,r)}S^y(t), \quad \dots \quad t_{m_L}=\argmax_{\overline{C(0,L)}\backslash \cup_{j=1}^{m_L-1}C(t_j,r)}S^y(t),
		\end{equation*}
		where $m_L$ is the number of steps until we cover~$\overline{C(0,L)}$. We denote the set  of candidate peaks by $T=\{t_1,\dots,t_{m_L}\}$.
		
		\item For $t\in T$ with observed height $S^y(t)$, compute the test values $p(S^y(t))$. The test value for identifying peaks is 
		\begin{equation}\label{eq: test value}
			\tilde{z}=M\cdot\max_{1\leq j\leq  M}(z\ast\ \tilde{\psi}_j)^2, \qquad p(u)=\mathbb{P}\left[ \max_{\overline{C\left(0,\frac{r}{2}\right)}}\tilde{z}>u\right],
		\end{equation} 
        where $M$ is the number of basis functions. 
		\item Apply a multiple testing procedure (Bonferroni~\cite{rupert2012simultaneous} or Benjamini-Hochberg~\cite{benjamini1995controlling}) on the set $\{p(S^y(t))\ | \ t\in T \}$, and declare as detected peaks those candidate peaks whose test values are higher than the threshold. Formally, the algorithm returns the set \begin{equation}\label{eq:Tu}
         T(u)=\left\{t\in T:\;S^y(t)\geq u\right\},   
        \end{equation}
        where $u$ is the threshold set by the multiple testing procedure.
	\end{enumerate}

\begin{figure}
	\captionsetup[subfigure]{labelformat=empty}
	\begin{subfigure}{.25\textwidth}
		\centering
		\includegraphics[height=3cm]{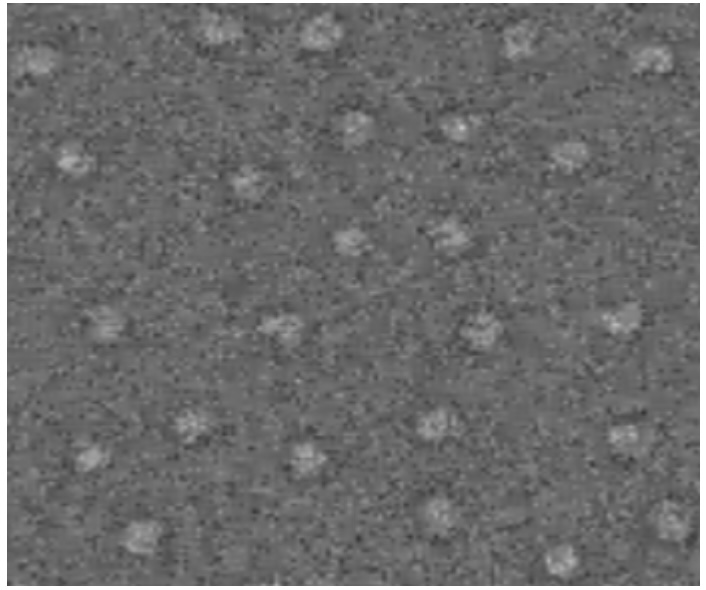}
		\caption{A}
	\end{subfigure}%
	\begin{subfigure}{.25\textwidth}
		\centering
		\includegraphics[height=3cm]{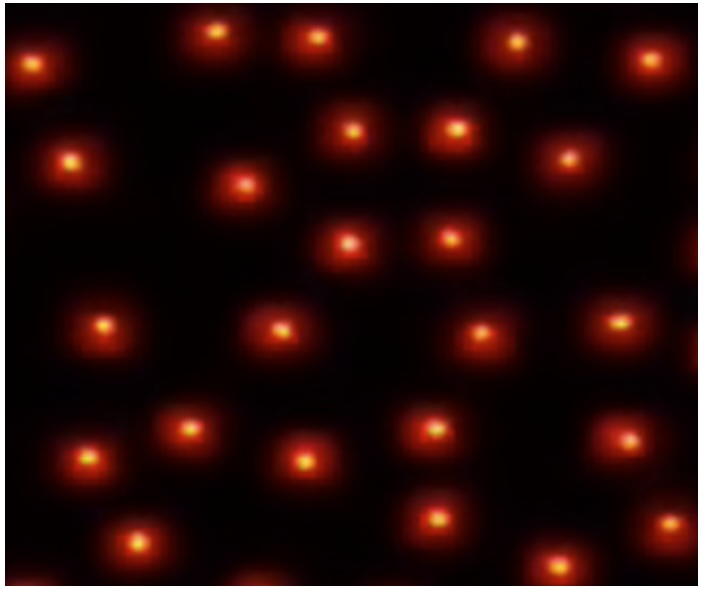}
		\caption{B}
	\end{subfigure}%
	\begin{subfigure}{.25\textwidth}
	\centering
	\includegraphics[height=3cm]{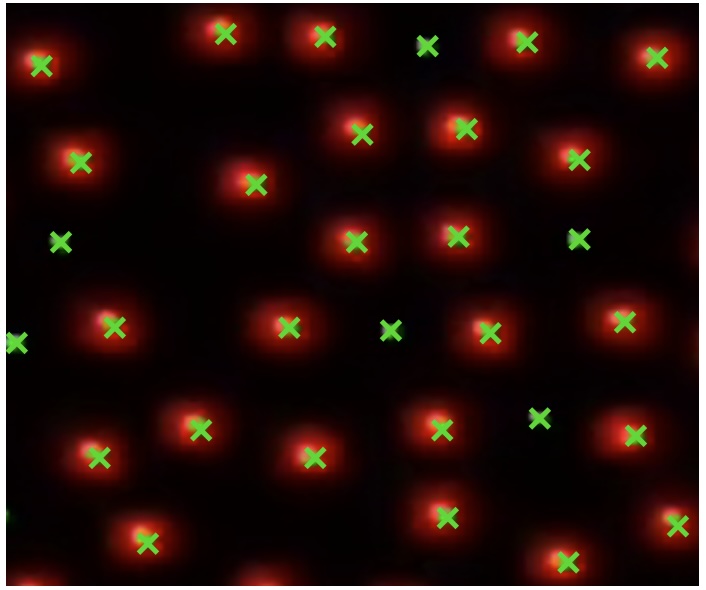}
	\caption{C}
\end{subfigure}%
\begin{subfigure}{.25\textwidth}
	\centering
	\includegraphics[height=3cm]{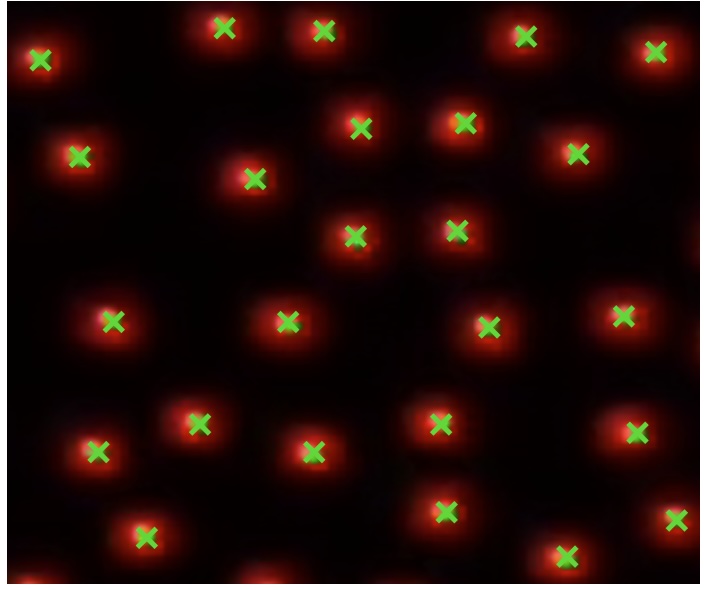}
	\caption{D}
\end{subfigure}%
\caption{Algorithm's outline (see Section \ref{sec:algorithm} for details). Given the noisy image A, the algorithm computes the scoring map B, selects candidate peaks (marked by green crosses  in C), and then returns the set of detected peaks (marked by green crosses  in D).}\label{fig:algorithm's flow}
\end{figure}

The intuition behind step~1 above is as follows. The expression~\eqref{eq: smoothed process} is the energy of the projection of the observed data onto the set $\left \{ \psi_{j} \right \}_{j=1}^{M}$ that spans the unknown objects. Since the noise is not spanned by this set of functions, we expect that projecting regions in the data containing objects would result in larger $S^{y}$ compared to areas that do not. We quantify this intuition in \Cref{appendix: higher SNR}. In light of step~1, in step $2$, higher values of the score map are expected to indicate the presence of an object.  We thus choose the highest peak in $S^{y}$ as the first candidate for an object location (and prove in \Cref{claim: the first N points in the algorithm are objects} that it indeed corresponds to an object with high probability). Since all pixels around this peak belong to that object, we remove them from subsequent peak searches, and repeat this process. The test value in step~3 allows to compute a detection threshold which distinguishes signal from noise by quantifying the expected peak height due to noise only. This is done in step~4 using a multiple hypothesis testing procedure.  We show in Section~\ref{sec:theorem} that the above procedure detects all objects, while controlling the error rate (under asymptotic conditions described in Section~\ref{sec:theorem}). The algorithm is illustrated in Figure~\ref{fig:algorithm's flow}.

\subsection{Theoretical guarantees}\label{sec:theorem}

To quantify the performance of the algorithm, we first need to rigorously define its power and error rate. As we would like to estimate the object centers~$\tau_{i}$ from~\eqref{eq:x}, we consider as a true positive (successful detection) a peak returned by our algorithm which is within distance~$\delta$ from~$\tau_{i}$, for some chosen parameter~$\delta$. We would obviously like~$\delta$ to be small, and discuss its choice below. We  define $ \mathbb{B}_1^{\delta} = \cup_{i=1}^NC(\tau_i,\delta)$ to be the union of the $\delta$~neighborhoods about the centers~$\tau_{i}$ (the true positive region), and the null region $ \mathbb{B}_0^{\delta}=\overline{C(0,L)}\backslash\mathbb{B}_1^\delta.$ Having applied the algorithm of Section~\ref{sec:algorithm}, the number of truly detected peaks and falsely detected peaks are
\begin{equation*}
		W(u)=\#\left\{T(u)\cap\mathbb{B}_1^{\delta}\right\}\quad \text{ and } \quad V(u)=\#\left\{T(u)\cap\mathbb{B}_0^\delta\right\},
\end{equation*}
where $T(u)$ is defined in~\eqref{eq:Tu}. Both  $W(u)$ and $V(u)$ are defined to be zero if $T(u)$ is empty.

The FWER is defined as the probability of obtaining at least one falsely detected peak, formally given by
\begin{align}\label{def: FWER}
\text{FWER}(u)=\mathbb{P}[V(u)\geq 1].
\end{align}
The false discovery rate ($\text{FDR}$) is defined as the expected proportion of falsely detected peaks, 
	\begin{equation}\label{def: fdr}
		\text{FDR}(u)= \mathbb{E}\left[\frac{V(u)}{\max\left\{W(u)+V(u),1\right\}}\right],
	\end{equation}
where the expectation is taken with respect to the noise. The power of the algorithm corresponds to the ratio of objects detected by the algorithm. Formally, for a given threshold $u>0$, we define
	\begin{align}\label{def: power}
		\text{Power}(u) &= \mathbb{E}\left[\frac{1}{N}\sum_{i=1}^{N}\mathbbm{1}_{\{T(u)\cap C(\tau_i,\delta)\neq\emptyset\}}\right]=\frac{1}{N}\sum_{i=1}^{N}\mathbb{P}[T(u)\cap C(\tau_i,\delta)\neq \emptyset].
	\end{align} 
	The indicator operator in \eqref{def: power} ensures that if more than one peak falls within the same $\delta$ neighborhood of some center, only one is counted, so that the power is not artificially inflated. The power is always non negative and less than $1$: power equal $0$ means that we did not detect any objects, while power equal $1$ means that we detected all objects. 

 The parameter $\delta$ that appears in the above definitions cannot be made arbitrarily small. A necessary condition for the algorithm to find the objects is that it detects the objects when no noise is present. This implies that  each peak detected in step~2 of the algorithm should be ``close'' to some $\tau_{i}$. The condition ensuring this is given in the following definition.
 \begin{definition}[Localization property]\label{def: localization property}	
    We say that $\delta > 0$ satisfies the localization property if, for every $1\leq i\leq N$, there exists some representative $q_i\in C(\tau_i,\delta)$ such that for all  $ t\in \overline{C(0,L)}\backslash \left(\cup_{j=1}^{i-1}C(\tau_j,2B) \cup_{j=i}^{N}C(\tau_j,\delta)\right)$
		\begin{equation*}
			S^x(q_i)-S^x(t)\geq \rho ||a_i||^2,
		\end{equation*}
		where $\rho$ is a positive constant, $a_{i}=\left ( a_{ij} \right )_{j=1}^{M}$ is the coefficients vector of $x_{i}$,  $a_{ij}$ is defined in~\eqref{eq:x}, and  $S^x = \sum_{j=1}^M \left(x\ast\tilde{\psi}_j\right)^2$ is the scoring map applied to the clean objects.
\end{definition}
We show in Appendix~\ref{appendix: delta zero estimation} how to compute a value of  $\delta$ which satisfies this property for a given set of functions $\{\psi_j\}_{j=1}^M$. Once~$\delta$ has been set, our goal is to compute an optimal detection threshold in step~4 of the algorithm that will control the error rate (FWER or FDR) below a pre-specified level~$\alpha$, while ensuring that the power is as close as possible to~$1$. We  make the following assumptions.

        \begin{assumption}\label{assumptions}  Let $\delta>0$ satisfy the localization property (\Cref{def: localization property}). Consider the following assumptions 
		\begin{enumerate}
            \item \label{assumption: distant centers} $\min_{1\leq l,j\leq N } \norm{\tau_l-\tau_j}_\infty>B+\frac{3}{2}\delta.$\label{assumption: distance between objectes}
	       \item $  \min_{t\in\partial\overline{C(0,L)}}\norm{\tau_j-t}_\infty \geq \frac{B}{2},\;\;1\leq j\leq N$, where $\partial \overline{C(t_0,L)}$ is the boundary of the hypercube.
			\item $\lim_{L\rightarrow\infty}\frac{N\cdot B^d}{L^d} = A$, where $ 0<A<1 $. \label{assumption: limit of the number of objects}
			\item Set $ a_L=\min_{1\leq i\leq N }\{||a_i||\}$, then $\frac{a_L^2}{L^d} \xrightarrow{L\rightarrow\infty} \infty$.
            \label{assumption: limit of snr}
		\end{enumerate}
	\end{assumption}
Assumption $1$ ensures that the objects are sufficiently separated (as a function of $\delta$). Assumption $2$ means that there are no objects near the boundary, as we want to detect only whole objects. Assumption $3$ states that the number of objects increases with $L$, while the density of the object occurrences over $\overline{C(0,L)}$ is fixed. Finally, assumption $4$ states that the norms of the objects are large enough to overcome the noise. 

We begin by controlling the FWER. Here, we set the threshold $u_{\rm{Bon}}$ to be the smallest $u$ such that $p(u)=\frac{\alpha}{M_L}$, where $M_L=\left\lceil\left(\frac{2L}{r}\right)^d\right\rceil$ (recall that $r$ is the parameter used in step~2 of the algorithm). Then, the following holds.

    \begin{theorem}\label{theorem:FWER} Suppose the assumptions above hold. Then, the algorithm of Section~\ref{sec:algorithm} with $r=2B+\delta$ and the Bonferroni threshold $u_{\rm{Bon}}$ satisfies	$$\lim_{L\rightarrow\infty}\text{Power}(u_{\rm{Bon}})=1, \quad \limsup_{L\rightarrow\infty}\text{FWER}(u_{\rm{Bon}})\leq\alpha.$$
	\end{theorem}
    \begin{proof}
        See Appendix~\ref{sec: proof FWER}.
    \end{proof}

To control the FDR below level~$\alpha$, we apply the Benjamini-Hochberg procedure \cite{benjamini1995controlling} in step~4 of the algorithm as follows:
	\begin{enumerate}
		\item Sort the set $p_{i}=\left\{p\left(S^y(t_i)\right)\big| t_i\in T\right\}$ in a non-decreasing order. 
		\item Find the largest index $k$ such that $p_{k}\leq\frac{k\alpha}{M_L} $, where $M_L=\left\lceil\left(\frac{2L}{r}\right)^d\right\rceil$. 
		\item Define the  threshold $u_{\rm{BH}}$ to be the smallest number $ u $ such that $p(u) = \frac{k\alpha}{M_L}$.
		\item Declare all $ k$ points in $T \left(u_{\rm{BH}}\right)$	as the detected objects' centers.
	\end{enumerate}
Since $u_{\rm{BH}}$ is a random variable, the expectation in~\eqref{def: fdr} and~\eqref{def: power}  is taken over all possible realizations of the random threshold $u_{\rm{BH}}$. We then have the following theorem.\       
     \begin{theorem}\label{theorem: fdr}
    Suppose the assumptions above hold.
    Then, the algorithm of Section~\ref{sec:algorithm} with $r=2B+\delta$  and the Benjamini-Hochberg threshold $u_{\rm{BH}}$ satisfies
 $$\lim_{L\rightarrow\infty}\text{Power}(u_{\rm{BH}})=1, \quad \limsup_{L\rightarrow\infty}\text{FDR}(u_{\rm{BH}})\leq\alpha.$$
	\end{theorem}
     \begin{proof}
        See Appendix~\ref{sec: proof FDR}.
    \end{proof}

We comment that the model~\eqref{eq: the model} was formulated in the continuous domain ($t$ is continuous), whereas in applications
we inevitably work with discretized data (such as images given by a discrete array of pixels). Our
analysis suggests that extending \Cref{theorem:FWER} to the discrete case (when $t$ is given on a
discrete finite set) is straightforward. In \Cref{theorem: fdr}, it is  relatively straightforward to prove that the power tends to one, but seems less trivial to prove
the control of the FDR, as the current proof requires the continuity of some intermediate functions.

\section{Numerical results}\label{sec: numerical results}
  
 \subsection{Simulations}\label{subsec: simulation studied}
We first evaluate the performance of the proposed algorithm using simulations. In each experiment, we simulate a $2$D realizations of $y$ of size $1024\times 1024$ pixels (i.e. $L=1024$) which contain images (objects) of size $64\times 64$ pixels ($B=64$). The objects are positioned at randomly chosen, non-overlapping locations, such that the density $\frac{N\cdot B^d}{L^d} = 0.5$. Each object is formed as  a linear combinations of $50$~Fourier Bessel functions~\cite{zhao2013fourier},  forcing symmetry so that the objects are real. The coefficients of each object are drawn independently from a uniform distribution over $[-1,1]$, such that the norm is~$1$. We use the method of Appendix~\ref{appendix: delta zero estimation} to conclude that $\delta= 10$ pixels satisfies the localization property of Definition~\ref{def: localization property} and set $r=2B+\delta$. The noise (a realization of $z$) is a centered Gaussian process with covariance kernel $C(x,y) = \sigma^2 e^{-2\norm{x-y}^2}$ , where $x,y\in \mathbb{R}^2$. As a pre-processing step, we simulate only the pure noise $ z $ in an image of size $202\times 202$ pixels, independently $N_{\text{sim}}=10^5$ times. We denote each noise realization by $z_i$. The image size of $z_i$ is chosen such that convolving it with the basis functions will give an output with the correct size needed to compute the samples $m_i = \max_{\overline{C\left(0,\frac{r}{2}\right)}}\tilde{z}_i$ defined in~\eqref{eq: test value}. 

\begin{figure}
		\centering
		\footnotesize
		\begin{minipage}[c]{\textwidth}
			\centering
   			\includegraphics[height=0.14\textheight]{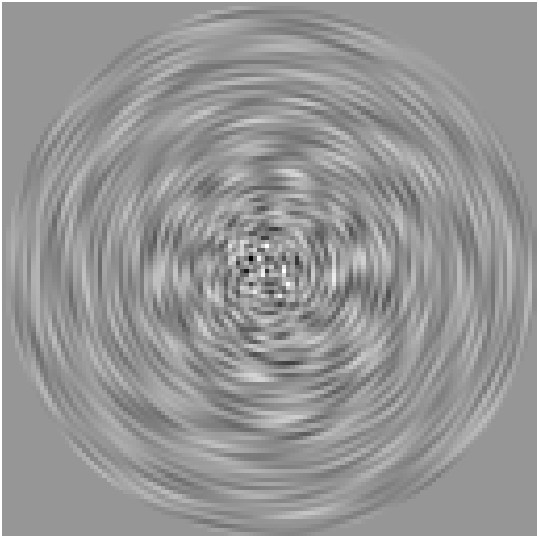}
			\includegraphics[height=0.14\textheight]{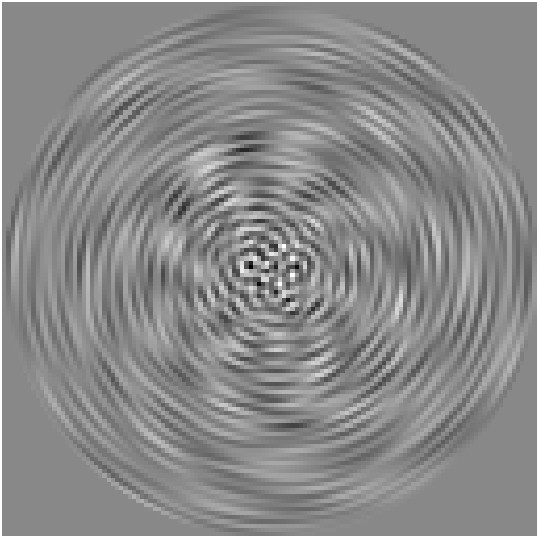}
			\includegraphics[height=0.14\textheight]{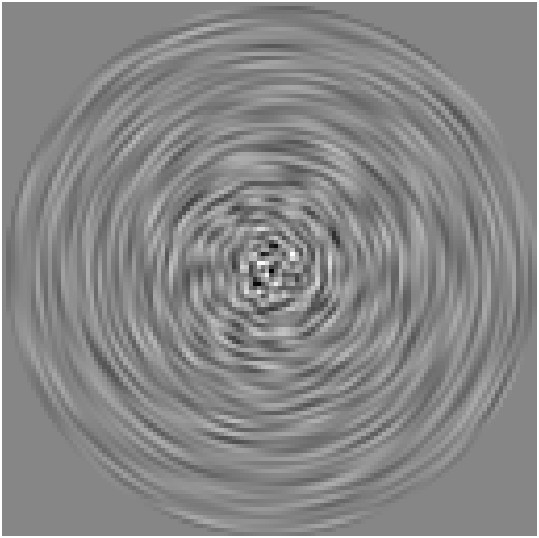}
			\includegraphics[height=0.14\textheight]{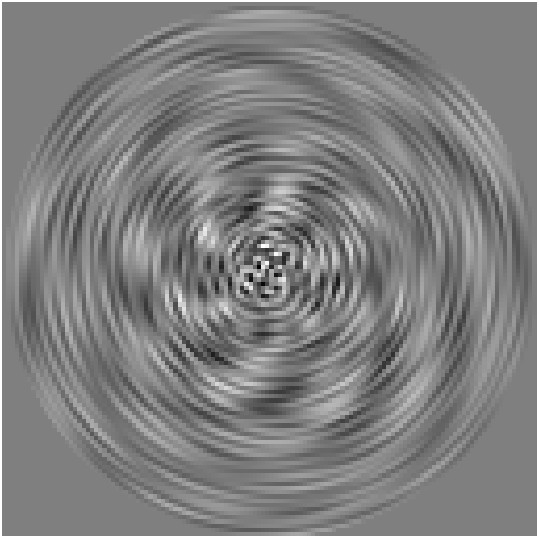}
			\vspace{0.5em}
		\end{minipage}
  		\begin{minipage}[c]{\textwidth}
			\centering
			\includegraphics[height=0.14\textheight]{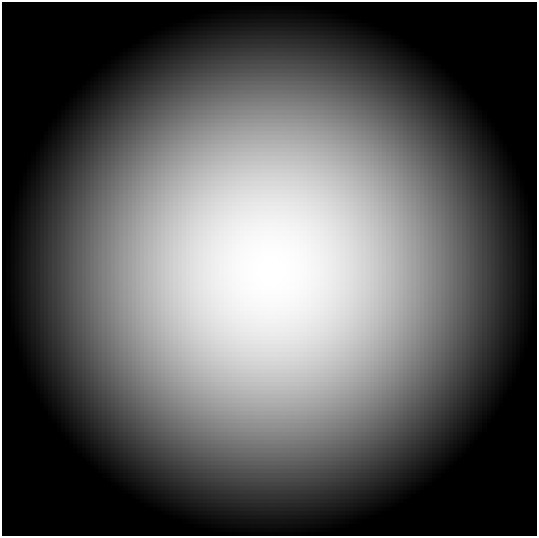}
			\includegraphics[height=0.14\textheight]{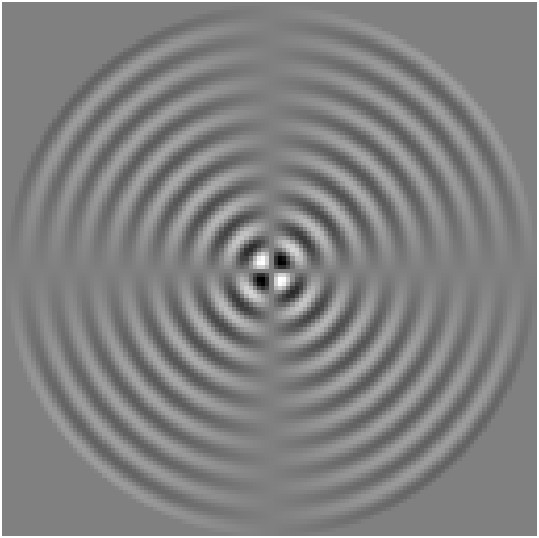}
			\includegraphics[height=0.14\textheight]{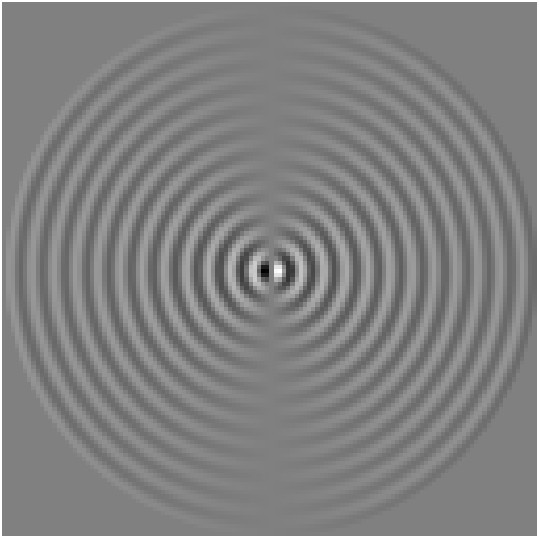}
			\includegraphics[height=0.14\textheight]{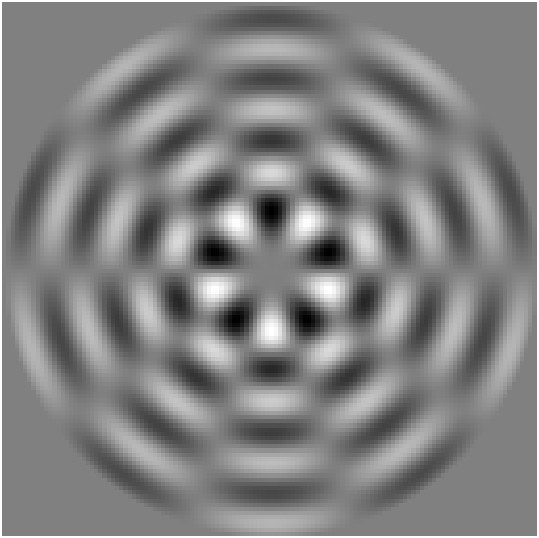}
			\vspace{0.5em}
		\end{minipage}
		\caption{Top row: four examples of objects, formed as linear combinations of~$50$ Fourier Bessel functions~\cite{zhao2013fourier},  with uniformly distributed coefficients, such that the norm of each object is~$1$. Bottom row: examples of the Fourier Bessel basis functions.}\label{fig: objects and basis simulations}
  \end{figure}

\begin{figure}
		\centering
		\footnotesize
		\begin{minipage}[c]{1\textwidth}
			\hspace{9em} Y
			\hspace{9.5em} Bonferroni 
			\hspace{7em} Benjamini-Hochberg 
		\end{minipage}
		\begin{minipage}[c]{\textwidth}
			\centering
			\rot{\hspace{4.5em}\vphantom{A}\vphantom{[}SNR = 1}
			\includegraphics[height=0.2\textheight]{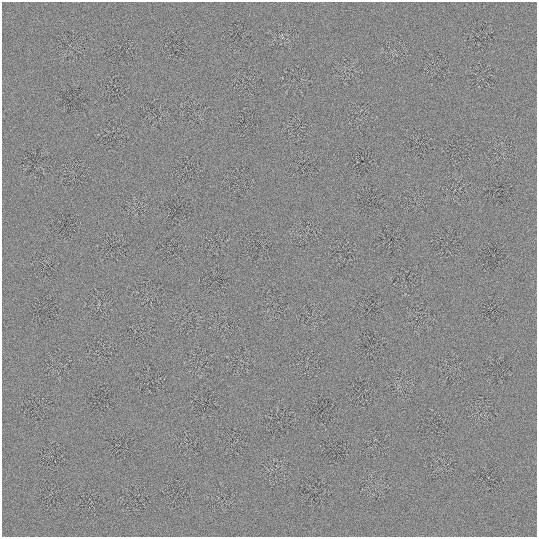}
			\includegraphics[height=0.2\textheight]{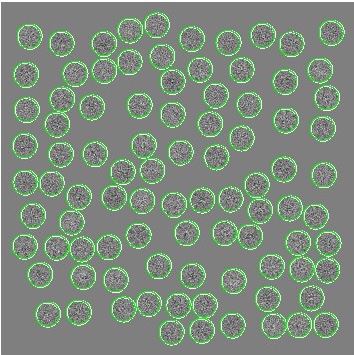}
			\includegraphics[height=0.2\textheight]{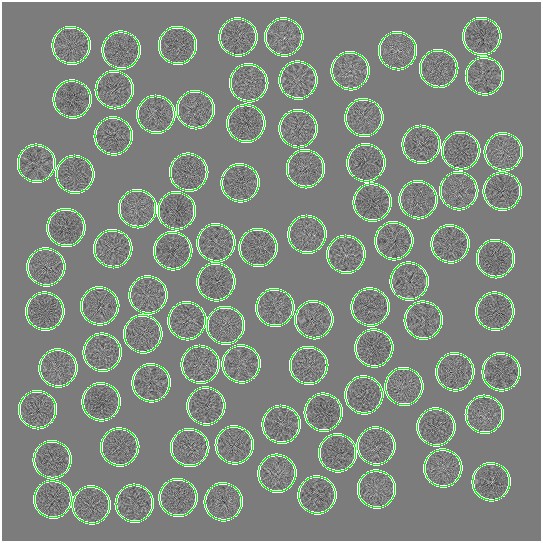}
			\vspace{0.5em}
		\end{minipage}
		\begin{minipage}[c]{\textwidth}
			\centering
		      \rot{\hspace{4em}\vphantom{A}\vphantom{[}SNR = 0.5}
			\includegraphics[height=0.2\textheight]{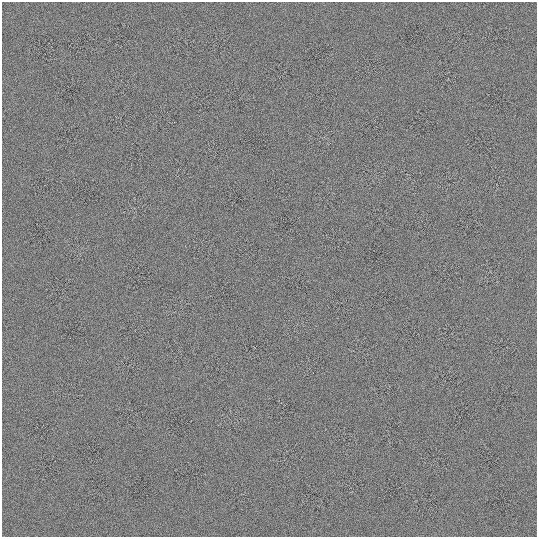}
			\includegraphics[height=0.2\textheight]{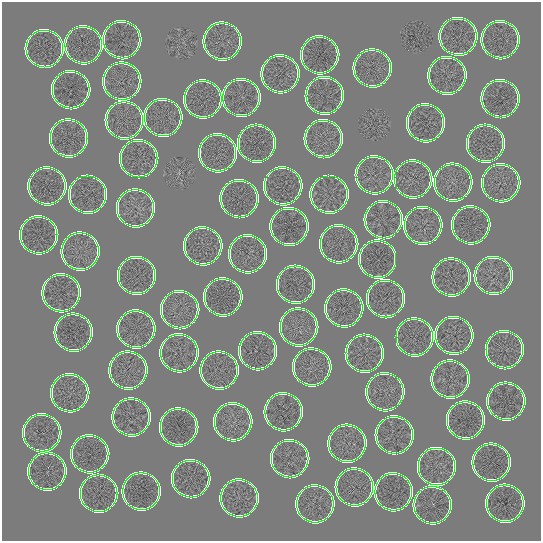}
			\includegraphics[height=0.2\textheight]{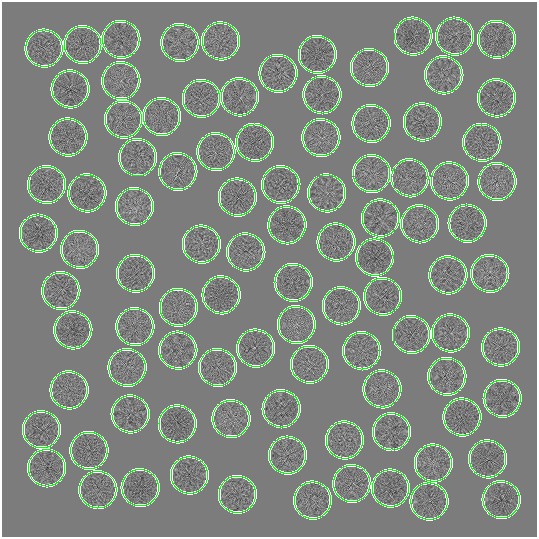}
			\vspace{0.5em}
		\end{minipage}
  		\begin{minipage}[c]{\textwidth}
			\centering
			\rot{\hspace{3.5em}\vphantom{A}\vphantom{[}SNR = 0.4}
			\includegraphics[height=0.2\textheight]{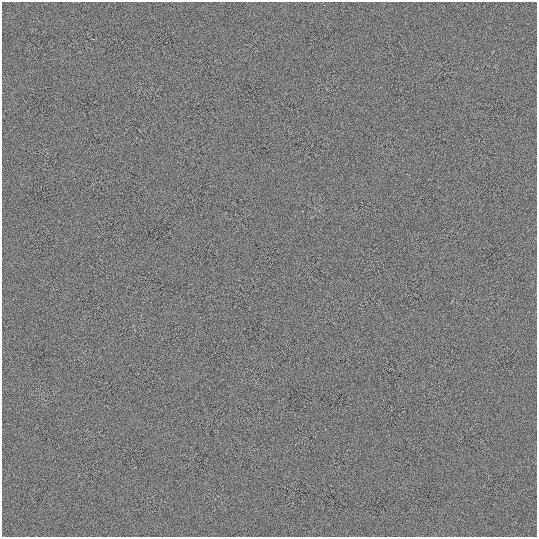}
			\includegraphics[height=0.2\textheight]{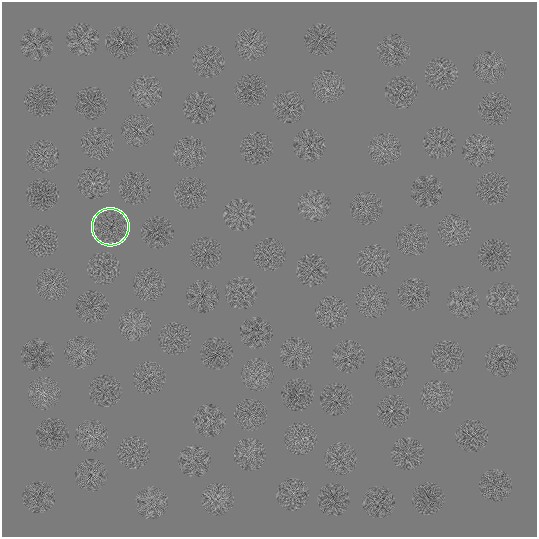}
			\includegraphics[height=0.2\textheight]{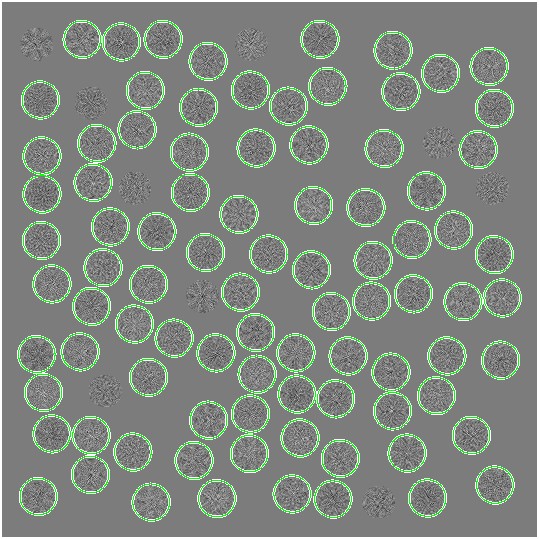}
			\vspace{0.5em}
		\end{minipage}
   		\begin{minipage}[c]{\textwidth}
			\centering
			\rot{\hspace{3.5em}\vphantom{A}\vphantom{[}SNR = 0.35}
			\includegraphics[height=0.2\textheight]{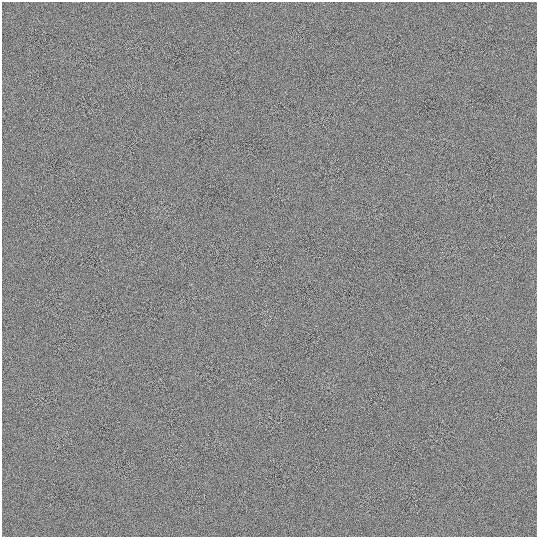}
			\includegraphics[height=0.2\textheight]{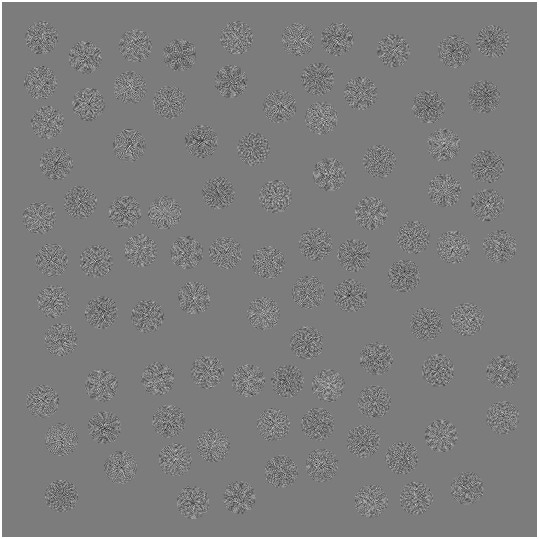}
			\includegraphics[height=0.2\textheight]{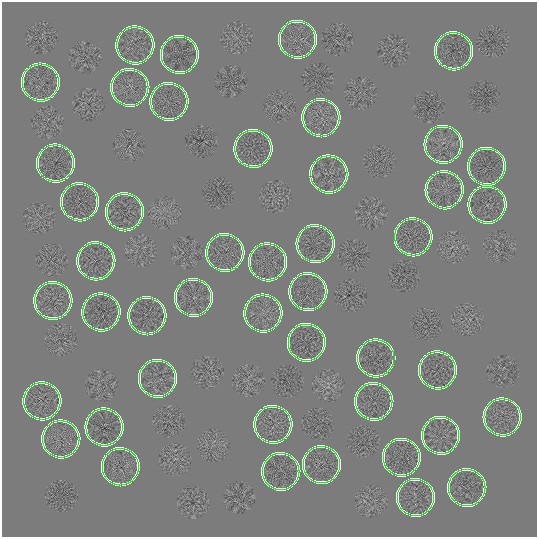}
			\vspace{0.5em}
		\end{minipage}
		\caption{Examples of detection results for different SNRs using the test value~\eqref{eq: test value}. The first column shows the noisy images $y$. The second and third columns show the results of our algorithm with respect to the Bonferroni and Benjamini-Hochberg procedures (see Theorems~\ref{theorem:FWER} and~\ref{theorem: fdr}). The detection results are shown on top of the noise-free images~$x$. The rows represent the SNRs at which the experiments have been conducted.}\label{fig: simulations examples}
	\end{figure}

From the realization of $ y $, we compute the candidate set $T$ by applying steps~$1$ and~$2$ of the algorithm (see Section~\ref{sec:algorithm}). 
  By the law of large numbers, using $p(u)=\mathbb{E}\left[\mathbbm{1}_{I(u)}\right]$, where $I(u)=\{\max_{\overline{C(0,\frac{r}{2}})}\tilde{z}>u\}$, we estimate $ p(u)$ defined in \eqref{eq: test value} as $\frac{\#\left\{i:m_i>u\right\}}{N_\text{sim}}$ , where $u=S^y(t)$ and $t\in T$. Using the estimates of $ p(u)$ we apply step $4$ of the algorithm (see Section~\ref{sec:algorithm}) with $\alpha = 0.05$. Finally, we estimate the FWER, FDR, and the power by repeating the entire experiment described in this paragraph $ 500 $ times independently, and using~\eqref{def: FWER}--\eqref{def: power}, where the expectation is estimated  by the sample mean computed over the 500 trials. To quantify the noise level in the data, we define the signal-to-noise-ratio (SNR)
 \begin{equation}\label{def: snr}
 \text{SNR}=\frac{\norm{a_L}^2}{\sigma^2B^{d}},  
 \end{equation} 
 where $a_L$ and $\sigma$ are defined in Assumption~\ref{assumption: limit of snr}. Note that since each of the objects has norm equal to~$1$, we have that~$\norm{a_L}=1$, and thus, $\sigma$ determines the $\text{SNR}$. Figure~\ref{fig: objects and basis simulations} shows examples of both the objects and basis functions. 
 Figure~\ref{fig: simulations examples} shows examples of the detection results for both the Bonferroni and Benjamini-Hochberg procedures, described in Theorems~\ref{theorem:FWER} and~\ref{theorem: fdr}.  Table~\ref{table: simulation results} presents the estimated FWER, FDR, and their power for different SNRs. As one can see, the error rates are always equal zero, indicating that the test value \eqref{eq: test value} is too conservative, which follows since it is given as~$M$ times the maximum of $\left\{\left(z\ast \psi_j\right)^2\right\}_{j=1}^{M} $. Instead, one can consider the more natural statistic $S^z = \sum_{j=1}^M\left(z\ast \psi_j\right)^2$, leading to an alternative test value
 \begin{equation}\label{eq: alternative test value}
 	\mathbb{P}\left[\max_{t\in\overline{C(0,\frac{r}{2}})}S^z(t)>u\right].
 \end{equation}
 Unfortunately, proving performance guarantees when using the test value~\eqref{eq: alternative test value} by our current proof strategy requires the continuity of~\eqref{eq: alternative test value} with respect to $u$, which we have not been able to establish.  Yet, we  demonstrate numerically the performance of the test value in~\eqref{eq: alternative test value}.  Table~\ref{table: simulation results altrnative test} presents the estimated FWER, FDR, and their power, for SNRs in which the test value defined in~\eqref{eq: test value} fails. Note that the error rates are not zero as in \Cref{table: simulation results} but still controlled.

\begin{table}
    \begin{center}
    \begin{tabular}{|c|c c|c c||} 
 \hline
 SNR & FWER & power  & FDR & power \\ [0.5ex] 
 \hline\hline
 $1$ & $0$ & $1$ & $0$ & $1$ \\ 
 \hline
 $0.5$ & $0$ & $0.98$ & $0$ & $1$ \\
 \hline
 $0.4$ & $0$ & $0.03$ & $0$ & $0.9$ \\
 \hline
 $0.35$ & $0$ & $0$ & $0$ & $0.15$ \\  
 \hline
\end{tabular}
\end{center}
\caption{FWER, FDR, and power, for different SNRs using the test value~\eqref{eq: test value}.} \label{table: simulation results}
\end{table}

\begin{table}
\begin{center}
\begin{tabular}{|c|c c|c c||} 
 \hline
 SNR & FWER & power  & FDR & power \\ [0.5ex] 
 \hline\hline
 $0.05$ & $0$ & $1$ & $0$ & $1$ \\
 \hline
 $0.03$ & $0.01$ & $0.95$ & $0.001$ & $0.99$ \\
 \hline
 $0.025$ & $0.02$ & $0.9$ & $0.002$ & $0.95$ \\
 \hline
 $0.02$ & $0.05$ & $0.66$ & $0.01$ & $0.88$ \\  
 \hline
\end{tabular}
\end{center}
\caption{FWER, FDR, and power for different SNRs using the test value~\eqref{eq: alternative test value}.} \label{table: simulation results altrnative test}
\end{table}

\subsection{Example of an experimental cryo-EM dataset}\label{subsec: data example}
\begin{figure}
		\centering
		\footnotesize
		\begin{minipage}[c]{\textwidth}
			\centering
   			\includegraphics[height=0.14\textheight]{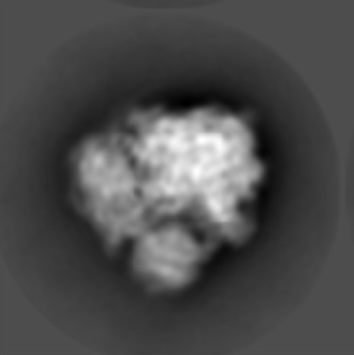}
			\includegraphics[height=0.14\textheight]{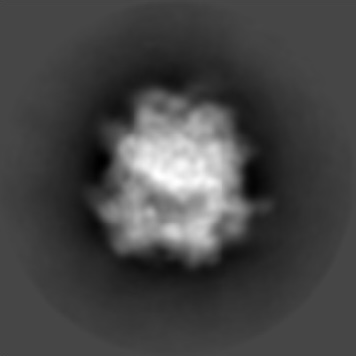}
			\includegraphics[height=0.14\textheight]{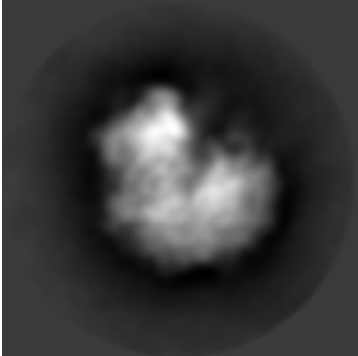}
			\includegraphics[height=0.14\textheight]{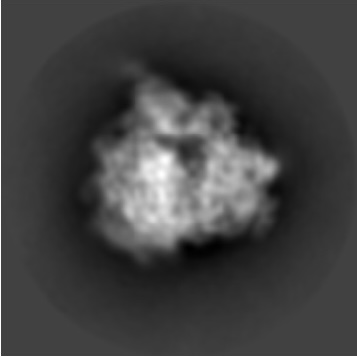}
			\vspace{0.5em}
		\end{minipage}
  		\begin{minipage}[c]{\textwidth}
			\centering
			\includegraphics[height=0.14\textheight]{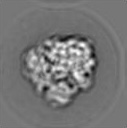}
			\includegraphics[height=0.14\textheight]{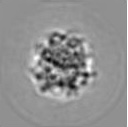}
			\includegraphics[height=0.14\textheight]{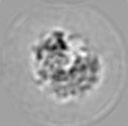}
			\includegraphics[height=0.14\textheight]{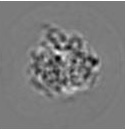}
			\vspace{0.5em}
		\end{minipage}
		\caption{Top row: examples of $4$ views generated by 2D classification~\cite{relion}. These are $4$ out of $34$ of the objects to detect.  Bottom row: examples of $4$ basis functions formed by Gram-Schmidt orthogonalization to the $34$ views.}\label{fig: objects and basis data example}
\end{figure}

We next demonstrate our algorithm for  template-based object detection in experimental data. In this setting, we are searching given templates in a large noisy image. The data consists of $161$ cryo-EM images of size $4096\times 4096$ pixels, taken from the Plasmodium Falciparum 80S ribosome (PF8r) data-set~\cite{10.7554/eLife.03080}, downloaded from the Electron Microscopy Public Image Archive~\cite{iudin2023empiar}. The size of the objects is approximately $360 \times 360$ pixels and there are between~$50$ to~$150$ objects in each image. Due to memory constraints, we down-sample the images such that the size of the objects after downsampling is $64\times 64$ pixels ($B=64 $). To generate templates, we used the algorithm of~\cite{relion}, resulting in $34$~different two-dimensional views of the underlying molecule. Four views are shown in Figure~\ref{fig: objects and basis data example}. The basis used by our algorithm is generated by applying the Gram-Schmidt orthogonalization to the~$34$ views. Four of the basis functions are shown in Figure~\ref{fig: objects and basis data example}. For the resulting basis, we set $\delta= 10$ pixels, and used the method of Appendix~\ref{appendix: delta zero estimation} to numerically verify that this~$\delta$ satisfies the localization property of Definition~\ref{def: localization property}. We also set $r=2B+\delta$. Next, we need to generate noise samples following the distribution of the noise in the data. To that end, we identified noise patches in the data, estimated their (isotropic) covariance matrix, and generated Gaussian samples with mean zero and the estimated covariance. All subsequent steps of the algorithm are implemented as in Section~\ref{subsec: simulation studied} with the test value \eqref{eq: alternative test value}.

\begin{figure}
		\centering
		\footnotesize
		\begin{minipage}[c]{\textwidth}
			\centering
   			\includegraphics[height=0.22\textheight]{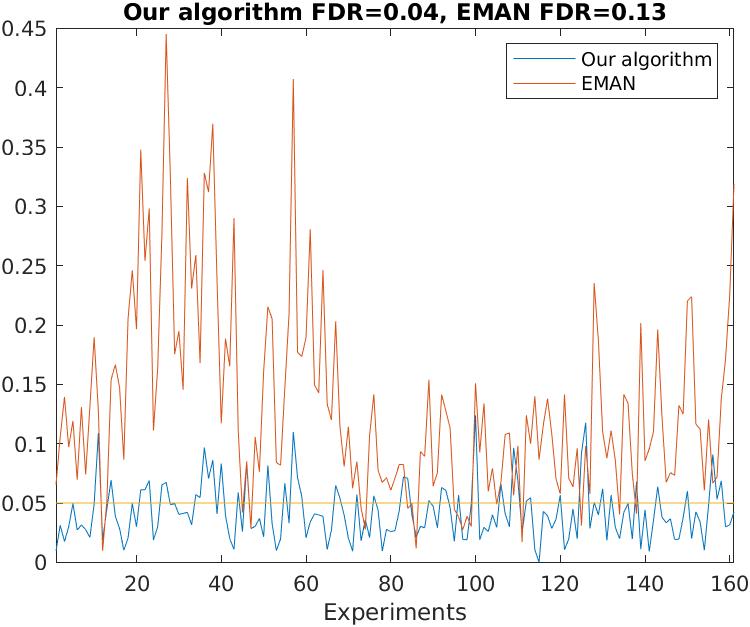}
			\includegraphics[height=0.22\textheight]{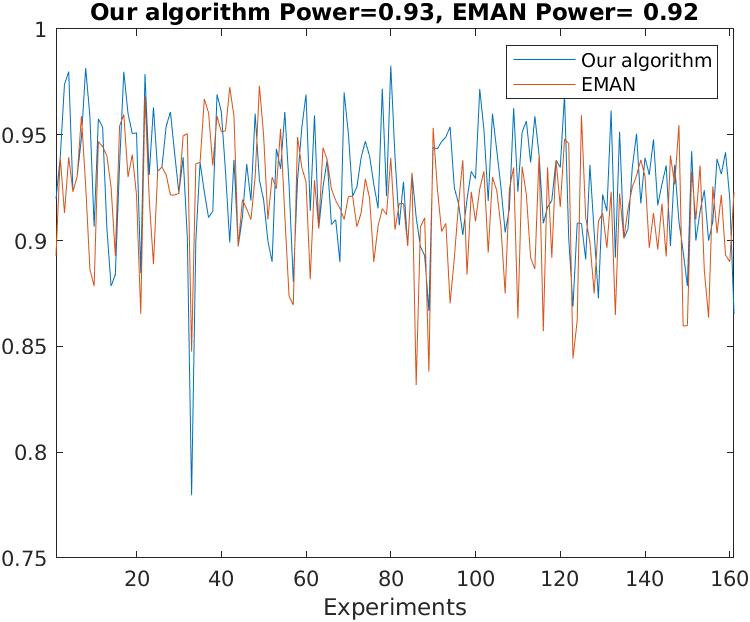}
			\vspace{0.5em}
		\end{minipage}
		\caption{Comparison of the FDR and power of our algorithm and the algorithm in EMAN. The left panel shows the FDP (i.e. $\frac{V(u)}{\max\left\{W(u)+V(u),1\right\}}$)  with $\alpha=0.05$ at the detection threshold per image (x-axis). The reported FDR is computed as the sample mean of FDP per images. Right panel shows the power at the detection threshold per image, i.e.  $\frac{1}{N}\sum_{i=1}^{N}\mathbbm{1}_{\{T(u)\cap C(\tau_i,\delta)\neq\emptyset\}}$ ($x$-axis). The reported power is computed as the sample mean of power per images.}\label{fig: fdr and power for real_data}
\end{figure}

\begin{figure}
		\centering
		\footnotesize
		\begin{minipage}[c]{\textwidth}
			\hspace{12em}EMAN
			\hspace{15.2em}Our algorithm
		\end{minipage}
		\begin{minipage}[c]{\textwidth}
			\centering
            \rot{\hspace{9em}\vphantom{A}\vphantom{[}93}
            \includegraphics[height=0.28\textheight]{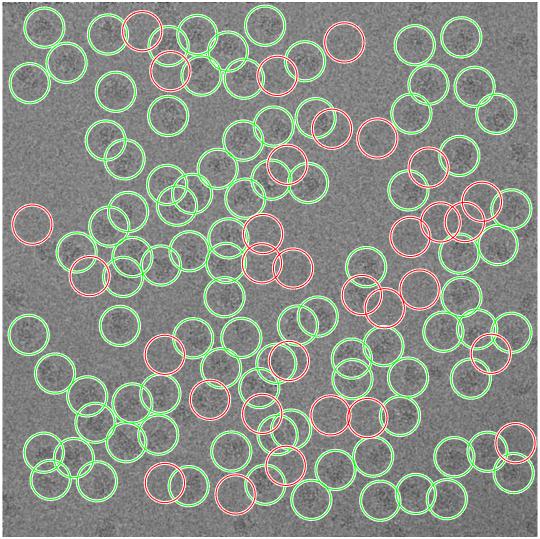}
            \includegraphics[height=0.28\textheight]{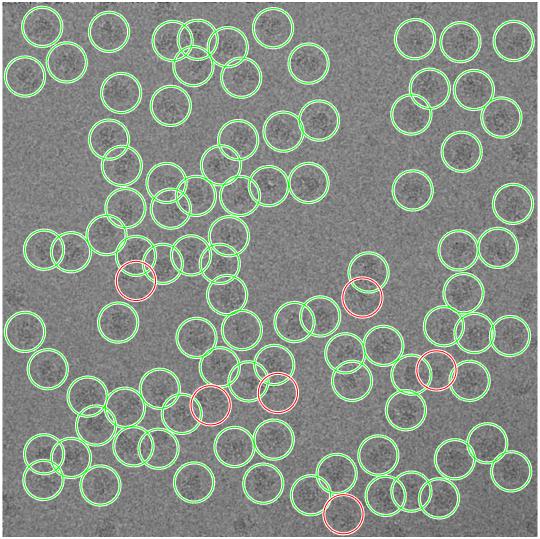}
		\end{minipage}
  
		\begin{minipage}[c]{\textwidth}
			\centering
            \rot{\hspace{9em}\vphantom{A}\vphantom{[}94}
            \includegraphics[height=0.28\textheight]{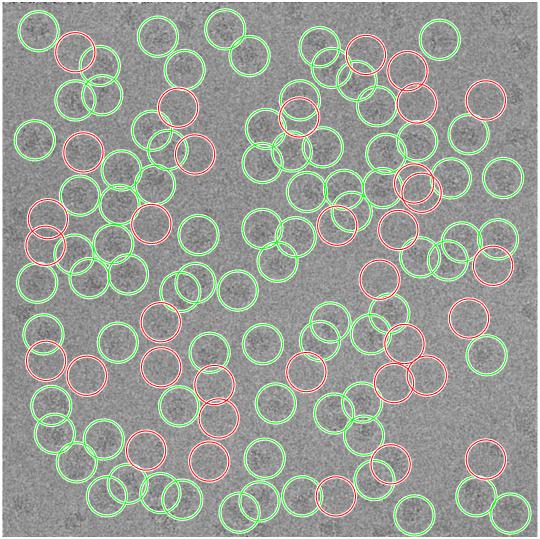}
            \includegraphics[height=0.28\textheight]{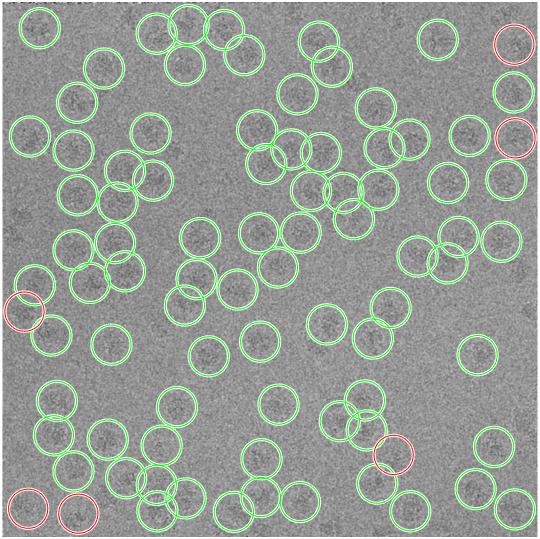}
		\end{minipage}
				\begin{minipage}[c]{\textwidth}
			\centering
            \rot{\hspace{9em}\vphantom{A}\vphantom{[}95}
            \includegraphics[height=0.28\textheight]{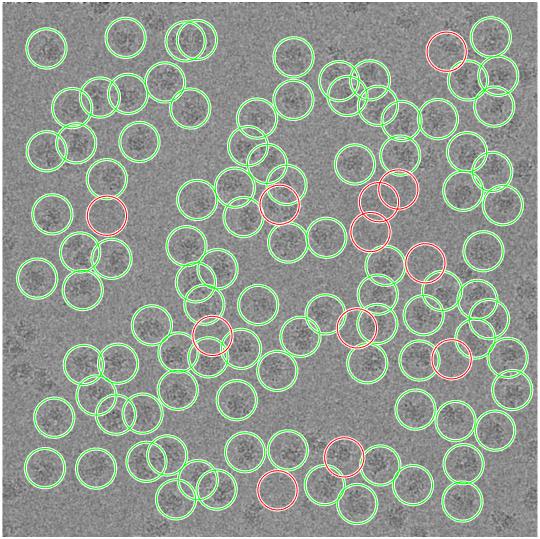}
            \includegraphics[height=0.28\textheight]{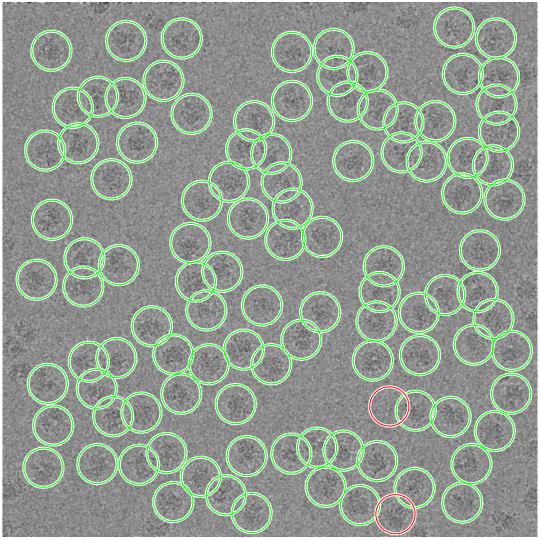}
		\end{minipage}
		
		\caption{Examples of detection results for the plasmodium falciparum 80S ribosome  dataset~\cite{10.7554/eLife.03080}. Each row corresponds to a different image with names provided by the benchmark data-set \cite{cryo-EMdataset}. The first and second columns display objects detected by EMAN and our algorithm, respectively. Green circles denote true positives, while red circles denote false positives.}
        \label{fig: comparison with EMAN}
	\end{figure}

We compare the results of our algorithm with the results of the particle picking algorithm in~EMAN~\cite{TANG200738} (one of the popular software packages for cryo-EM data processing). We use the same input to both algorithms. The detection threshold of the algorithm in EMAN is computed automatically by the software.  To evaluate the results of both algorithms we used the benchmark dataset~\cite{cryo-EMdataset}, which reports the centers of the objects in the data. Since the centers' identification process in~\cite{cryo-EMdataset} is manual, it is subject to small errors. Therefore, we define true positive whenever an algorithm detects a point in the data that is up to $90$ pixels from the center reported in the benchmark dataset (this value corresponds to 50\% of the radius of the object). Figure~\ref{fig: fdr and power for real_data} 
presents the FDR and power for both algorithms. As one can see, the FDR of our algorithm is controlled, lower than EMAN's and our power is higher. We do not present results for the FWER, as both algorithms make more than one mistake on each image (experiment),  which leads to an FWER of $1$. Finally, Figure~\ref{fig: comparison with EMAN} shows $3$ examples of the detection results. Note that our algorithm has higher number of true positives (green circles) than EMAN, while having fewer false positives (red circles).


\section{Discussion and future research}\label{sec:future}

In this work, we presented an object detection algorithm that allows to control the power and error rates simultaneously. The performance bounds of the algorithm are proven in an asymptotic regime, though the numerical experiments show that it performs well in real-world scenarios, and in particular, in the non-asymptotic regime. 

An important future work is deriving non-asymptotic performance bounds. In particular, we conjecture that we can replace the requirement that the norms of the objects tend to infinity by requiring that the projected SNR (see~\ref{appendix: higher SNR}) is sufficiently large.  In addition, as noted in Section~\ref{subsec: simulation studied}, the test value we use seems to be sub-optimal; we intend to analyze the performance of the proposed algorithm with the alternative test value~\eqref{eq: alternative test value}. The main difficulty is proving that this test value is continuous with respect to~$u$. Finally, we wish to design bases that allow~$\delta$ to be as small as possible. 
\section*{Acknowledgments}
Amitay Eldar was supported by NSF-BSF award 2019733.
Keren Mor Waknin was supported in part by a fellowship from the Edmond J. Safra Center for Bioinformatics at Tel-Aviv University and The Jacob and Riva Damm Endowment Fund. Tamir Bendory  supported in part by BSF grant no. 2020159, NSF-BSF grant no. 2019752, and  ISF grant no.1924/21. Samuel Davenport and Armin Schwartzman were supported in part by NIH/NIBIB Award R01EB026859. Yoel Shkolnisky was supported in part by NIH/NIGMS Award R01GM136780-01.

\begin{appendices}

  \section{Auxiliary claims}
    \begin{claim}\label{claim: centers norms} Let $l\in\{1,\dots,N\}$ and $S^x = \sum_{j=1}^M \left(x\ast\tilde{\psi}_j\right)^2$, then,
		$ S^x(\tau_l) = ||a_l||^2$.
	\end{claim}	
    \begin{proof}
		Denote $ \psi^i_j(t)=\psi_j(t-\tau_i) $, then,
		\begin{align*}
			S^x(\tau_l) &= \sum_{j=1}^M \left(x\ast\tilde{\psi}_j\left(\tau_l\right)\right)^2 =\sum_{j=1}^M \left(\sum_{i=1}^Nx_i\ast\tilde{\psi}_j\left(\tau_l\right)\right)^2\\\nonumber
			&=\sum_{j=1}^M \left(\sum_{i=1}^N\sum_{k=1}^Ma_{ik}\psi^i_k \ast\tilde{\psi}_j\left(\tau_l\right)\right)^2\\\nonumber
			&=\sum_{j=1}^M \left(\sum_{i=1}^N\sum_{k=1}^Ma_{ik}\int_{\overline{C(0,L)}}\psi^i_k(s)\tilde{\psi}_j\left(\tau_l-s\right)ds\right)^2 \\\nonumber
			&=\sum_{j=1}^M \left(\sum_{i=1}^N\sum_{k=1}^Ma_{ik}\int_{\overline{C(0,L)}}\psi_k(s-\tau_i)\psi_j\left(s-\tau_l\right)ds\right)^2\\\nonumber
			&=\sum_{j=1}^M \left(\sum_{i=1}^N\sum_{k=1}^Ma_{ik}\int_{\overline{C(0,L)}}\psi_k(t)\psi_j\left(t-\left(\tau_l-\tau_i\right)\right)dt\right)^2\\\nonumber
			&=\sum_{j=1}^M \left(\sum_{i=1}^N\sum_{k=1}^Ma_{ik}\int_{C(0,B)}\psi_k(t)\psi_j\left(t-\left(\tau_l-\tau_i\right)\right)dt\right)^2\\\nonumber
			&=\sum_{j=1}^M \left(\sum_{i=1}^N\sum_{k=1}^Ma_{ik}\int_{C(0,B)}\psi_k(t)\psi_j\left(t\right)\delta_{il}dt\right)^2\\\nonumber
			&=\sum_{j=1}^M \left(\sum_{i=1}^N\sum_{k=1}^Ma_{ik}\delta_{kj}\delta_{il}\right)^2=\sum_{j=1}^M a_{lj}^2.
		\end{align*}
  where the fifth equality is due to change of variables and the seventh equality follows from Assumption \ref{assumption: distance between objectes}.
	\end{proof}

 \begin{claim}[Projected SNR is higher at an object's center]\label{appendix: higher SNR} Let $l\in\{1,\dots,N\}$. Denote by $\text{SNR}_{y}(\tau_{l})=\frac{\norm{a_l}^2}{\sigma^2B^{d}}$ the input signal-to-noise ratio and by $\text{SNR}_{S^y}(\tau_{l}) = \frac{S^x(\tau_{l})}{\mathbb{E}\left[S^z(\tau_{l})\right]}$ the projected signal-to-noise ratio, where $S^z = \sum_{j=1}^M \left(z\ast\tilde{\psi}_j\right)^2$ . Then, $\text{SNR}_{S^y}(\tau_l) = \gamma \times \text{SNR}_{y}(\tau_{l})$ for some $\gamma\ge 1$.  
 \end{claim} 
 \begin{proof}
 	We can write $ z\vert_{C(0,B)} = \sum_{j=1}^\infty b_j\psi_j$ where $ \left\{\psi_j\right\}_{j=M+1}^\infty$ is an orthonormal complement of $\{\psi_j\}_{j=1}^{M}$ and $ b_j $ are random variables.  Denote  $\sigma^2=\mathbb{E}\left[z^2(t)\right]$. Then,
	\begin{equation*}
	    \sigma^2B^d = \mathbb{E}\left[\int_{C(0,B)}z^2(t)dt\right]= \mathbb{E}\left[\sum_{j=1}^{\infty} b_j^2\right]=\sum_{j=1}^{\infty}\mathbb{E}\left[b_j^2\right],
	\end{equation*}
    where the second equality follows from Parseval's theorem, and the third follows from the monotone convergence theorem.
    On the other hand, by stationarity, for each $1 \leq l \leq N$, 
		\begin{align*}
			\mathbb{E}\left[S^z(\tau_l)\right]&=\mathbb{E}\left[S^z(0)\right] = \mathbb{E}\left[\sum_{j=1}^M \left(z\ast\tilde{\psi}_j\right)^2(0)\right]\\
			&=\mathbb{E}\left[\sum_{j=1}^M \left(\int_{C(0,B)}z(u)\psi_j(u)du\right)^2\right]\\
			&=\sum_{j=1}^M \mathbb{E}\left[b^2_j\right] =\frac{1}{\gamma} \sigma^2B\leq \sigma^2B,
		\end{align*} 
  where $\gamma = \frac{\sum_{j=1}^M \mathbb{E}\left[b^2_j\right]}{\sum_{j=1}^\infty \mathbb{E}\left[b^2_j\right]} $.
	The latter, combined with \Cref{claim: centers norms}, implies that $$ \text{SNR}_{S^y}(\tau_l)  = \gamma\frac{S^x(\tau_{l})}{\sigma^2B^d} = \gamma\frac{\norm{a_l}^2}{ \sigma^2B^d}= \gamma\text{SNR}_{y}(\tau_l).$$
 \end{proof}
 \begin{claim}\label{appendix: mixed term definition}
     $  S^y = S^z+S^x + H$, where   $S^z = \sum_{j=1}^M \left(z\ast\tilde{\psi}_j\right)^2$, $S^x = \sum_{j=1}^M \left(x\ast\tilde{\psi}_j\right)^2$, $H = 2\sum_{j=1}^{M}\left(x\ast\tilde{\psi}_j\right)\cdot\left(z\ast\tilde{\psi}_j\right)$ .
 \end{claim}
 \begin{proof}
          \begin{align*}
		S^y &= \sum_{j=1}^{M}\left(y\ast\tilde{\psi}_j\right)^2=\sum_{j=1}^{M}\left((x+z)\ast\tilde{\psi}_j\right)^2\\\nonumber
		&=\sum_{j=1}^{M}\left(z\ast\tilde{\psi}_j\right)^2 + \sum_{j=1}^{M}\left(x\ast\tilde{\psi}_j\right)^2 +2\sum_{j=1}^{M}\left(x\ast\tilde{\psi}_j\right)\cdot\left(z\ast\tilde{\psi}_j\right)\\\nonumber
		&=S^z+S^x + H.
	\end{align*}
 \end{proof}
	\begin{claim}[Controlling the growth of the mixed term $H$]\label{claim: bound for W}
		Let $t\in \mathbb{R}^d$, and denote by $I(t)=\{i \mid ||\tau_i-t||_\infty<B\}$ the set of indices of centers $\tau_i$ such that $||\tau_i - t||_\infty< B$. Then, $$\left|H(t)\right|\leq c\cdot\max_{i\in I(t)}\{||a_i||\}\cdot\sqrt{S^z(t)},$$
		where   $S^z = \sum_{j=1}^M \left(z\ast\tilde{\psi}_j\right)^2$, $H = 2\sum_{j=1}^{M}\left(x\ast\tilde{\psi}_j\right)\cdot\left(z\ast\tilde{\psi}_j\right)$ and $c$ is a positive constant.
	\end{claim}
 \begin{proof}
		Let $t\in\mathbb{R}^d$, then, as in \Cref{claim: centers norms}, by a change of variable we get
		\begin{align*}
			H(t)&=2\sum_{j=1}^{M}\left(\sum_{i=1}^{N}\sum_{k=1}^{M}a_{ik}\psi^i_k\ast\tilde{\psi}_j\right)\left(z\ast\tilde{\psi}_j\right)\left(t\right)\\
			&=2\sum_{j=1}^M \left(\sum_{i=1}^N\sum_{k=1}^Ma_{ik}\int_{C(0,B)}\psi_k(s)\psi_j\left(s-\left(t-\tau_i\right)\right)ds\right)\left(z\ast\tilde{\psi}_j\right)\left(t\right).\\\nonumber
		\end{align*}
        Hence, 
		\begin{align*}
			H(t)&=2\sum_{j=1}^M \left(\sum_{i\in I(t)}\sum_{k=1}^Ma_{ik}\int_{C(0,B)}\psi_k(s)\psi_j\left(s-\left(t-\tau_i\right)\right)ds\right)\left(z\ast\tilde{\psi}_j\right)\left(t\right).
		\end{align*}
		By the Cauchy–Schwartz inequality, 
		\begin{align*}
			|H(t)|\leq 2\sum_{j=1}^M\sum_{i\in I(t)}\sqrt{\sum_{k=1}^{M}a_{ik}^2}\sqrt{\sum_{k=1}^M\left(\int_{C(0,B)}\psi_k(s)\psi_j\left(s-\left(t-\tau_i\right)\right)ds\right)^2}\cdot\left|\left(z\ast\tilde{\psi}_j\right)\left(t\right)\right|.
		\end{align*}
		Applying the Cauchy–Schwartz inequality again, and using the fact that $||\psi_k||=1$, for  all $1 \leq k \leq M$, we get that $$\int_{C(0,B)}\left|\psi_k(s)\psi_j\left(s-\left(t-\tau_i\right)\right)\right|ds\leq 1.$$
        Note that for any $t\in\mathbb{R}^d$, since the centers are at least at distance $B$ from each other, there can be at most $3^d-1$ centers satisfying $||t-\tau_i||_{\infty}<B$, so $|I(t)|\leq 3^d-1$. Hence,
		\begin{align*}
			 \left|H(t)\right|&\leq 2\sum_{j=1}^{M}(3^d-1)\cdot\max_{i\in I(t)}\{||a_i||\}\cdot \sqrt{M}|(z\ast\tilde{\psi}_j)|(t)\\&
			=2\cdot (3^d-1)\max_{i\in I(t)}\{||a_i||\}\cdot\sqrt{M}\sum_{j=1}^{M}|(z\ast\tilde{\psi}_j)|(t)
			\\&\leq2M\cdot (3^d-1)\max_{i\in I(t)}\{||a_i||\}\sqrt{\sum_{j=1}^{M}(z\ast\tilde{\psi}_j)^2(t)}\\&=c\cdot\max_{i\in I(t)}\{||a_i||\}\cdot\sqrt{S^z(t)}.
		\end{align*}
	\end{proof}
	\begin{definition}\label{def: lesssim}
		Let $a, b, x, y \in \mathbb{R}$. We say $ae^b \lesssim xe^y$ if there exist positive constants $c_1$ and $c_2$ that do not depend on $L$, such that $ae^b \leq c_1 xe^{c_2 y}$.
	\end{definition}
    \begin{claim}[Controlling the growth of the projected noise]
    \label{claim: tail prob. estimate}
		For every $u>0$, we have $$\mathbb{P}\left[\max_{\overline{C(0,L)}} S^z>u\right]\leq\mathbb{P}\left[\max_{\overline{C(0,L)}} \tilde{z}>u\right] \lesssim L^de^{-u},$$
     	where $S^z = \sum_{j=1}^M \left(z\ast\tilde{\psi}_j\right)^2$,  $\tilde{z}=M\cdot\max_{1\leq j\leq M}\left(\max_{\overline{C(0,L)}}\left(z\ast\tilde{\psi}_j\right)^2\right)$ and $\lesssim$ is defined in \ref{def: lesssim}.
	\end{claim}
    \begin{proof}
            Denote for all $1\leq j\leq M$, $z_j=z\ast \tilde{\psi}_j$. The first inequality follows immediately by \begin{align}
            	 S^z(t)=\sum_{j=1}^Mz_j^2(t)\leq M\cdot\max_{1\leq j\leq M}\left(\max_{\overline{C(0,L)}}z_j^2(t)\right)=\tilde{z}(t)
            \end{align}  surely for all $t\in\mathbb{R}^d$. Then,
	\begin{align*}
		\mathbb{P}\left[\max_{\overline{C(0,L)}} S^z>u\right]&=\mathbb{P}\left[\max_{\overline{C(0,L)}}\sum_{j=1}^{M}z_j^2>u\right]\leq\mathbb{P}\left[M\cdot\max_{1\leq j\leq M}\left(\max_{\overline{C(0,L)}}z_j^2\right)>u\right]\\&=\mathbb{P}\left[\max_{\overline{C(0,L)}}\tilde{z}>u\right]
        \leq \mathbb{P}\left[\bigcup_{j=1}^{M}\left(\max_{\overline{C(0,L)}}z_j^2\right)>\frac{u}{M}\right]\\&\leq \sum_{j=1}^{M}\mathbb{P}\left[\left(\max_{\overline{C(0,L)}}z_j^2\right)>\frac{u}{M}\right],
	\end{align*}
	where the last inequality is the union bound.
	Let us bound the probability $$\mathbb{P}\left[\max_{\overline{C(0,L)}}z_j^2>\frac{u}{M}\right].$$
    To this end, divide $\overline{C(0,L)}$ into $L^d$ boxes of volume 1. Then,
    \begin{align*}  \mathbb{P}\left[\max_{\overline{C(0,L)}}z_j^2>\frac{u}{M}\right]\leq L^d\mathbb{P}\left[\max_{\overline{C(0,1)}}z_j^2>\frac{u}{M}\right]&=L^d\mathbb{P}\left[\max_{\overline{C(0,1)}}|z_j|>\sqrt{\frac{u}{M}}\right].
    \end{align*}
    Now, set $$\mu_j=\mathbb{E}\left[\max_{\overline{C(0,1)}}z_j\right], \quad \sigma_j^2=\Var\left[z_j\right].$$
    By the Borell–TIS inequality \cite{Adler2007} we get,
    \begin{align*}  \mathbb{P}\left[\max_{\overline{C(0,1)}}|z_j|>\sqrt{\frac{u}{M}}\right]\leq 2e^{-\frac{\left(\sqrt{\frac{u}{M}}-\mu_j\right)^2}{2\sigma_j^2}},
    \end{align*}
    so that 
    \begin{align*}  \mathbb{P}\left[\max_{\overline{C(0,L)}}z_j^2>\frac{u}{M}\right]\leq 2L^de^{-\frac{\left(\sqrt{\frac{u}{M}}-\mu_j\right)^2}{2\sigma_j^2}}\leq 2L^d e^{\frac{-\left(\frac{u}{M}+\eta^2\right)+2\sqrt{\frac{u}{M}}\mu}{2\sigma'^2}},
    \end{align*}
    where $\mu=\max_{1\leq j\leq M}\mu_j$, $\eta=\min_{1\leq j\leq M}\mu_j$, and $\sigma'^2=\max_{1\leq j\leq M}\sigma_j^2$. 
    Overall, we get 
    \begin{align*}  \mathbb{P}\left[\max_{\overline{C(0,L)}}\tilde{z}>\frac{u}{M}\right]\leq 2ML^de^{\frac{-\left(\frac{u}{M}+\eta^2\right)+2\sqrt{\frac{u}{M}}\mu}{2\sigma'^2}} \lesssim L^de^{-u}.
    \end{align*} 
    \end{proof}

    \begin{claim}\label{claim: contiuity of p}
        The test value defined in \eqref{eq: test value} i.e. $p(u)=\mathbb{P}\left[\max_{\overline{C(0,\frac{r}{2}})}\tilde{z}>u\right]$  is continuous for all $u\in\mathbb{R}$.
    \end{claim}
    \begin{proof}
        By \cite{ylvisaker1968note}, for each $1\leq j\leq M$, $z_j=z\ast\tilde{\psi}_j$ is a continuous random variable. For all $1\leq j\leq M$, 
        \begin{align*}
            \max_{\overline{C(0,\frac{r}{2}})}z_j^2=            \left(\max_{\overline{C(0,\frac{r}{2}})}\left|z_j\right|\right)^2,
        \end{align*} which implies the continuity of each of the $            \max_{\overline{C(0,\frac{r}{2}})}z_j^2$. 
        Then, 
        \begin{align*}
            \max_{\overline{C(0,\frac{r}{2}})}\tilde{z}=\max_{\overline{C(0,\frac{r}{2}})}\left(M\cdot\max_{1\leq j\leq M}z_j^2\right)=M\cdot\max_{1\leq j\leq M}\max_{\overline{C(0,\frac{r}{2}})}z_j^2
        \end{align*}
        which yields the continuity of $            \max_{\overline{C(0,\frac{r}{2}})}\tilde{z}$ as the maxima of $M$ continuous random variables. 
        
    \end{proof}
    
	\begin{claim}[All objects are picked first]\label{claim: the first N points in the algorithm are objects} Let the assumptions from  \Cref{assumptions} hold, and suppose that our algorithm is applied with $r=2B+\delta$. Denote by $E_L^{(1)}$ the event that the first $ N $ points in $ T $ are in $ \mathbb{B}_{1}^{\delta}$. Then
	\begin{align}
		\mathbb{P}\left[\left(E_L^{(1)}\right)^\complement\right] \lesssim 2^N\cdot L^de^{-a_L^2}\xrightarrow{L\rightarrow\infty}0
	\end{align}
	where $\lesssim$ is defined in \ref{def: lesssim}.
	\end{claim}
	\begin{proof}
		Recall that $T=\{t_i\}_{i=1}^{m_L}$, and
		denote $A_L^{(i)}=\{t_i\in\mathbb{B}_1^{\delta}\}$. Note that $E_L^{(1)}=\bigcap\limits_{i=1}^{N}A^{(i)}_L$.  We will prove that for all $1\leq i \leq N$, $$\mathbb{P}\left[\left(A_L^{(i)}\right)^{\complement}\right]\lesssim 2^{i-1}\cdot L^de^{-a_L^2}.$$ This will complete the proof 
		as $$\mathbb{P}\left[\left(\bigcap\limits_{i=1}^{N}{A}^{(i)}_L\right)^{\complement}\right]\leq\sum_{i=1}^{N}\mathbb{P}\left[\left(A_L^{(i)}\right)^{\complement}\right]\lesssim 2^N\cdot L^de^{-a_L^2}\xrightarrow{L\rightarrow\infty}0, $$ where the limit tends to zero by Assumptions \ref{assumption: limit of the number of objects} and \ref{assumption: limit of snr}.
        Assume (without loss of generality) that $||a_1||\geq||a_2||\geq \dots \geq ||a_N||$.
		We start with $\left(A_L^{(1)}\right)^{\complement}=\{t_1\in \mathbb{B}_0^{\delta}\}$. Recall that $q_1\in C(\tau_1,\delta)$ is the representative defined in Definition \ref{def: localization property}, and that $t_1=\argmax_{t\in \overline{C(0,L)}}{S^y(t)}$. Note that
		\begin{align*}
			\{t_1\in \mathbb{B}_0^{\delta}\}&=\{S^y(t_1)-S^y(q_
			1)\geq0,t_1\in\mathbb{B}_0^{\delta}\}\\\nonumber
			&=\{S^z(t_1)+S^x(t_1)+H(t_1)-S^z(q_1)-S^x(q_1)-H(q_1)\geq0,t_1\in \mathbb{B}_0^{\delta}\}\\
			&\subseteq \{2\max_{t\in\overline{C(0,L)}}S^z(t)+|H(t_1)|+|H(q_1)|\geq S^x(q_1)-S^x(t_1), t_1\in\mathbb{B}_0^{\delta}\}.
		\end{align*}
		Now, using   \Cref{claim: bound for W}, we get that 
		\begin{align*}
			|H(t_1)|\lesssim \max_{i\in I(t_1)}\{||a_i||\}\sqrt{S^z(t_1)} ,\quad |H(q_1)|\lesssim \max_{i\in I(q_1)}\{||a_i||\}\sqrt{S^z(q_1)}.
		\end{align*}
		Altogether we get, 
		\begin{align*}
			|H(t_1)|+|H(q_1)|&\lesssim ||a_1||\left(\sqrt{S^z(t_1)}+\sqrt{S^z(q_1)}\right)\\&\lesssim ||a_1||\sqrt{\max_{t\in\overline{C(0,L)}}S^z(t)}.
		\end{align*}
		The latter yields, 
		\begin{align*}
			\{t_1\in\mathbb{B}_0^{\delta}\}&\subseteq \{\max_{t\in\overline{C(0,L)}}S^z(t)+||a_1||\sqrt{\max_{t\in\overline{C(0,L)}}S^z(t)}\gtrsim S^x(q_1)-S^x(t_1), t_1\in\mathbb{B}_0^{\delta}\}\\&\subseteq \{\max_{t\in\overline{C(0,L)}}S^z(t)+||a_1||\sqrt{\max_{t\in\overline{C(0,L)}}S^z(t)}\gtrsim \rho||a_1||^2\},
		\end{align*}
		where we used the localization property (\Cref{def: localization property}), since $t_1\in\mathbb{B}_0^{\delta}$.
		This implies        
            \begin{align*}
            \mathbb{P}\left[t_1\in\mathbb{B}_0^{\delta}\right]&\leq\mathbb{P}\left[\max_{\overline{C(0,L)}}S^z(t)+ ||a_1||\sqrt{\max_{t\in\overline{C(0,L)}}S^z(t)} \gtrsim \rho||a_1||^2\right]\\&  \leq\mathbb{P}\left[\max_{\overline{C(0,L)}}S^z(t)\gtrsim \rho||a_1||^2 \right]+\mathbb{P}\left[\sqrt{\max_{t\in\overline{C(0,L)}}S^z(t)} \gtrsim \rho||a_1||\right].
		\end{align*}
		Using \Cref{claim: tail prob. estimate}, we get that for all $u>0$, $$\mathbb{P}\left[\max_{\overline{C(0,L)}}S^z(t)>u\right]\lesssim L^de^{-u},$$ 
		which yields,  $$\mathbb{P}\left[\left(A_L^{(1)}\right)^{\complement}\right]\lesssim L^de^{-||a_1||^2}.$$ 
		For ${A_L^{(2)}}^{\complement}$, we are going to intersect it with the event $A_L^{(1)}$ which means we want to bound the probability that $t_2\in\mathbb{B}_0^{\delta}$ and $t_1\in\mathbb{B}_1^{\delta}$. Note that $t_1\in\mathbb{B}_1^{\delta} \iff t_1\in C(\tau_i,\delta)$ for a unique  $i\in\{1,\dots,N\}$.
		Denote $||a_j||=\max_{n\neq i}\{||a_n||\}$. Since $t_2=\argmax_{\overline{C(0,L)}\backslash C(t_1,r)}S^y,$ 
		then, on the event $A_L^{(1)}$, we have that $q_j\in \overline{C(0,L)}\backslash C(t_1,r)$, which implies
		\begin{align*}
			\{{A_L^{(2)}}^{\complement}\cap A_L^{(1)}\}&=\{t_2\in \mathbb{B}_0^{\delta},t_1\in \mathbb{B}_1^{\delta}\}\\&= \{S^y(t_2)-S^y(q_
			j)\geq0,t_2\in\mathbb{B}_0^{\delta},t_1\in C(\tau_i,\delta)\}\\\nonumber
			&=\{S^z(t_2)+S^x(t_2)+H(t_2)-S^z(q_j)-S^x(q_j)-H(q_j)\geq0,t_2\in\mathbb{B}_0^{\delta},t_1\in C(\tau_i,\delta)\}\\
			&\subseteq \{2\max_{t\in\overline{C(0,L)}}S^z(t)+|H(t_2)|+|H(q_j)|\geq 	S^x(q_j)-S^x(t_2), t_2\in\mathbb{B}_0^{\delta},t_1\in C(\tau_i,\delta)\}.
		\end{align*}
		Now, in order to bound $|H(t_2)|$ and $|H(q_j)|$ with $||a_j||$, first observe that $I(q_j)=\{j\}$. Indeed, since $||\tau_j-q_j||_\infty<\frac{\delta}{2}$ and $||\tau_i-\tau_j||_\infty>B+\frac{3}{2}\delta$ for all $i,j$, it holds that $I(q_j)=\{j\}$.   Second, for $I(t_2)$, note that since $t_2\notin C(t_1,r)$ and $ t_1\in C(\tau_i,\delta)$, then, $||t_2-\tau_i||_\infty\geq \left|\ ||t_2-t_1||_\infty-||t_1-\tau_i||_\infty \ \right|>\frac{r}{2}-\frac{\delta}{2}=B+\frac{\delta}{2}-\frac{\delta}{2}=B$, which implies that $i\notin I(t_2)$. Hence $\max_{n\in I(t_2)}\{||a_n||\}\leq ||a_j||$. Combining the latter with  \Cref{claim: bound for W}, we get
		\begin{align*}
			|H(t_2)|+|H(q_j)|&\lesssim ||a_j||\left(\sqrt{S^z(t_2)}+\sqrt{S^z(q_j)}\right)\\&\lesssim ||a_j||\sqrt{\max_{t\in\overline{C(0,L)}}S^z(t)}.
		\end{align*}
		In order to bound $S^x(q_j)-S^x(t_2)$ we split to cases. If $\tau_i=\tau_1$, then, $||a_j||=\max_{n\neq 1}||a_n||=||a_2||$ and $t_2\notin C(\tau_1,2B)\cup_{j=1}^NC(\tau_j,\delta)$, so it follows from \Cref{def: localization property} that $$S^x(q_2)-S^x(t_2)\geq \rho||a_2||^2.$$ 
	Otherwise, $\tau_i\neq \tau_1$, which means that $||a_j||=||a_1||$ and $t_2\notin C(\tau_i,2B)\cup_{j=1}^NC(\tau_j,\delta)$, which implies that $t_2\notin\cup_{j=1}^NC(\tau_j,\delta)$, and hence, we can apply \Cref{def: localization property}  again to get $$S^x(q_1)-S^x(t_2)\geq \rho||a_1||^2.$$ To conclude, in both cases we have $S^x(q_j)-S^x(t_2)\geq \rho||a_j||^2$, and thus, 
		\begin{align*} \mathbb{P}\left[\left(A_L^{(2)}\right)^{\complement}\cap A_L^{(1)}\right]&\leq \mathbb{P}\left[\max_{\overline{C(0,L)}}S^z(t)+||a_j||\sqrt{\max_{t\in\overline{C(0,L)}}S^z(t)}\gtrsim \rho||a_j||^2\right]\\&\lesssim L^de^{-||a_j||^2}.
		\end{align*}
	
		Overall, we get
		\begin{align*}
			\mathbb{P}\left[\left(A_L^{(2)}\right)^{\complement}\right]&=  \mathbb{P}\left[\left(A_L^{(2)}\right)^{\complement}\cap A_L^{(1)}\right]+ \mathbb{P}\left[\left(A_L^{(2)}\right)^{\complement}\cap \left(A_L^{(1)}\right)^{\complement}\right]\\&
			\lesssim L^de^{-||a_j||^2} +\mathbb{P}\left[\left(A_L^{(1)}\right)^{\complement}\right]
			\\&
			 \lesssim L^de^{-||a_j||^2}+L^de^{-||a_k||^2}\lesssim 2L^de^{-a_L^2}.
		\end{align*}
		Then, we can repeat this argument for all $3\leq i\leq N$ and induction to get 
	\begin{align*}
		\mathbb{P}\left[\left(A_L^{(i)}\right)^{\complement}\right]&=  \mathbb{P}\left[\left(A_L^{(i)}\right)^{\complement}\cap \left(\cap_{j=1}^{i-1} A_L^{(j)}\right)\right]+ \mathbb{P}\left[\left(A_L^{(i)}\right)^{\complement}\cap \left(\cap_{j=1}^{i-1} A_L^{(j)}\right)^{\complement}\right]\\& \lesssim L^de^{-a_L^2}+\sum_{j=1}^{i-1} \mathbb{P}\left[\left(A_L^{(j)}\right)^\complement\right]\\& \lesssim L^de^{-a_L^2} + L^de^{-a_L^2}\sum_{j=1}^{i-1}2^{j-1}\\&= L^de^{-a_L^2} + L^de^{-a_L^2}(2^{i-1}-1)=2^{i-1}L^de^{-a_L^2}.
	\end{align*}
	\end{proof}
	\begin{claim}[The multiple testing procedure stops after more than or equal to $N$ steps]\label{claim: BH algorithm have more then N steps}
 Let the assumptions from \Cref{assumptions} hold, and suppose that our algorithm is applied with  $r=2B+\delta$. Denote by $E_L^{\text{Bon}}$, $E_L^{\text{BH}}$ the events that the Bonferroni or Benjamini-Hochberg procedures, defined in \Cref{theorem:FWER,theorem: fdr}, stop  after more than or equal to $N$ steps. Then, 
	\begin{align}
	\mathbb{P}\left[\left(E_L^{\text{Bon}}\right)^\complement\right], \mathbb{P}\left[\left(E_L^{\text{BH}}\right)^\complement\right] \lesssim 2^N\cdot L^de^{-a_L^2} + N\cdot\frac{1}{L^{4d}}\xrightarrow{L\rightarrow\infty}0
	\end{align}
	where $\lesssim$ is defined in \ref{def: lesssim}.
	\end{claim}
	\begin{proof}
		We start with the Benjamini-Hochberg procedure. 
		Recall that $k$ is the largest index such that $$p\left(S^y(t_k)\right)\leq\frac{k}{M_L}\alpha.$$
		Let $E_L^{\text{BH}}=\{\text{The BH procedure stops after more or equal than N steps}\}$. Recall that $E_L^{(1)}$ is the event that the first $ N $ points in $ T $ are in $ \mathbb{B}_{1}^{\delta}$. In \Cref{claim: the first N points in the algorithm are objects} we proved that $\mathbb{P}\left[\left(E_L^{(1)}\right)^{\complement}\right]\lesssim 2^N\cdot L^de^{-a_L^2} $. Set~$ \tilde{E_L}^{(1)}=\cap_{i=1}^N\{S^z(q_i)\leq \log(L^{4d})\}$, where $q_i$ defined in \ref{def: localization property}. By \Cref{claim: tail prob. estimate}, and the union bound it holds that $\mathbb{P}\left[\left(\tilde{E_L}^{(1)}\right)^{\complement}\right]\lesssim N\cdot\frac{1}{L^{4d}}$. Combining both we have that  $$\mathbb{P}\left[\left(E_L^{(1)}\right)^{\complement}\cup \left(\tilde{E_L}^{(1)}\right)^{\complement}\right]\lesssim 2^N\cdot L^de^{-a_L^2} + N\cdot\frac{1}{L^{4d}}\xrightarrow{L\to \infty} 0,$$ 
		where the limit is zero due to Assumptions \ref{assumption: limit of the number of objects} and \ref{assumption: limit of snr}. 
		We shall prove that $\exists L_0>0$ such that for all $L>L_0$, $$E_L^{(1)}\cap \tilde{E_L}^{(1)}\subset E_L^{\text{BH}},$$ which finishes the proof. Assume without loss of generality that $||a_1||\geq||a_2||\geq \dots \geq ||a_N||$. Let $\omega\in E_L^{(1)}\cap \tilde{E_L}^{(1)}$. Then, $S^y(t_N)(\omega)\geq S^y(q_N)(\omega)$ (see definition of $t_N$ in Section \ref{sec:algorithm}).
		Now,
		\begin{align*}
			p(S^y(t_N)(\omega))&=\mathbb{P}\left[\max_{\overline{C\left(0,\frac{r}{2}\right)}}\tilde{z}>S^y(t_N)(\omega)\right]\leq\mathbb{P}\left[\max_{\overline{C\left(0,\frac{r}{2}\right)}}\tilde{z}>S^y(q_N)(\omega)\right]\\
			&= \mathbb{P}\left[\max_{\overline{C\left(0,\frac{r}{2}\right)}}\tilde{z}>S^x(q_N)+S^z(q_N)(\omega)+H(q_N)(\omega)\right]
			\\&\leq 
			\mathbb{P}\left[\max_{\overline{C\left(0,\frac{r}{2}\right)}}\tilde{z}>S^x(q_N)-|H(q_N)(\omega)|\right]
				\\&\leq \mathbb{P}\left[\max_{\overline{C\left(0,\frac{r}{2}\right)}}\tilde{z}>\rho||a_N||^2-\max_{i\in I(q_N)}\{||a_i||\}\sqrt{S^z(q_N)(\omega)}\right]
				\\&= \mathbb{P}\left[\max_{\overline{C\left(0,\frac{r}{2}\right)}}\tilde{z}>\rho||a_N||^2-||a_N||\sqrt{S^z(q_N)(\omega)}\right]
			\\&\leq \mathbb{P}\left[\max_{\overline{C\left(0,\frac{r}{2}\right)}}\tilde{z}\gtrsim\rho||a_N||^2(1-\frac{\sqrt{\log(L^d)}}{||a_N||}\right]
			\\&\lesssim \left(\frac{r}{2}\right)^de^{-a_L^2}\overset{\forall L>L_0}{\leq} \frac{N}{M_L}\alpha,
		\end{align*} 
        where the third equality follows from \Cref{appendix: mixed term definition},  the forth inequality is due to non-negativity of $S^z$ and the fifth inequality follows by \Cref{claim: bound for W} and \Cref{def: localization property}. The six equality is due to the fact that $I(q_N)=\{N\}$, the seven inequality follows from the definition of  $\tilde{E_L}^{(1)}$. In the last two inequalities, we used \Cref{claim: tail prob. estimate} and Assumptions \eqref{assumption: limit of the number of objects},\eqref{assumption: limit of snr}. To conclude we proved that for every $L>L_0$ and $\omega\in E_L^{(1)}\cap \tilde{E_L}^{(1)}$, the BH procedure stops  after more than or equal to $N$ steps which means by definition that $\omega\in E_L^{\text{BH}}$.
        
		For the Bonferroni procedure, we use the last inequality to get $$p(S^y(t_N)(\omega))\lesssim \left(\frac{r}{2}\right)^de^{-a^2_L}\overset{\forall L>\tilde{L}_0}{\leq}\frac{1}{M_L}\alpha,$$ 
		which implies that $\exists \tilde{L}_0$ such that $\forall L>\tilde{L}_0$ it holds that $E_L^{(1)}\cap \tilde{E_L}^{(1)}\subset E_L^{\text{Bon}}$.
	\end{proof}

     Recall that $V(u)=\#\left\{T(u)\cap\mathbb{B}_0^{\delta} \right\}$ for all $u\in\mathbb{R}$, where $T(u)$ and $\mathbb{B}_0^{\delta}$ are defined in \Cref{sec:algorithm,sec:theorem}. Define for all $u\in\mathbb{R}$   
     \begin{align}\label{definition: U(u)}
     	 U(u)=\sum_{j=1}^{M_L}\mathbbm{1}_{I_j(u)}
     \end{align} where $M_L=\left\lceil\left(\frac{2L}{r}\right)^d\right\rceil$, $I_j(u)=\left\{\max_{\overline{C\left(v_j,\frac{r}{2}\right)}} \tilde{z}>u\right\}$, and $\overline{C(0,L)}=\cup_{j=1}^{M_L}\overline{C(v_j,\frac{r}{2})}$ is a partition of $\overline{C(0,L)}$ into evenly sized closed hypercubes of side length $\frac{r}{2}$ (assuming without loss of generality that $L$ is divided by $\frac{r}{2}$). 
    \begin{claim}\label{claim: V less U} Let the assumptions from section \eqref{assumptions} hold, and suppose that  our algorithm is applied with $r=2B+\delta$.
    On the event $E_L^{(1)}$(defined in Claim \ref{claim: the first N points in the algorithm are objects}) where the $N$ first points in~$T$ are in~$\mathbb{B}_1^{\delta}$, it holds that
	$$V(u)\leq U(u),$$
 for all $u>0$. 
\end{claim}
\begin{proof}
	Note that $ V(u)=\#\left\{T(u)\cap\mathbb{B}_0^{\delta} \right\} \leq \sum_{j=1}^{M_L}\#\left\{\overline{C\left(v_j,\frac{r}{2}\right)}\cap T(u)\cap\mathbb{B}_0^{\delta}\right\} $. The algorithm in \Cref{sec:algorithm} gives that if $ t_i,t_j\in  T$,  then $ \norm{t_i-t_j}_\infty \geq \frac{r}{2}$, and therefore, $\forall v_j\in \mathbb{R}^d $ there is only one point from $ T $ that can be in $ \overline{C(v_j,\frac{r}{2})} $. The latter yields that $\#\left\{T(u)\cap\mathbb{B}_0^{\delta} \right\} \leq 1$ and thus $$ V(u) \leq \sum_{j=1}^{M_L}\mathbbm{1}_{\left\{\exists t\in \overline{C(v_j,\frac{r}{2})}\cap T(u)\cap\mathbb{B}_0^{\delta}\right\}}= \sum_{j=1}^{M_L}\mathbbm{1}_{\left\{\exists t\in \overline{C(v_j,\frac{r}{2})}\cap T\cap\mathbb{B}_0^{\delta},S^y(t)>u\right\}}.$$\\
	On the event $E_L^{(1)}$, we have that $T\cap \mathbb{B}_0^{\delta}=\{t_i\}_{i=N+1}^{m_L}$, then, again by the algorithm in \Cref{sec:algorithm}, for all $t\in T\cap \mathbb{B}_0^{\delta}$:   $S^y(t) =  S^z(t)$.
	Hence, on the event $E_L^{(1)}$, 
\begin{align*}
 V(u)&\leq\sum_{j=1}^{M_L}\mathbbm{1}_{\left\{\exists t\in \overline{C(v_j,\frac{r}{2})}\cap T\cap\mathbb{B}_0^{\delta},S^z(t)>u\right\}}\\&\leq \sum_{j=1}^{M_L}\mathbbm{1}_{\left\{\max_{\overline{C\left(v_j,\frac{r}{2}\right)}} S^z(t)>u\right\}}\\&\leq\sum_{j=1}^{M_L}\mathbbm{1}_{\left\{\max_{\overline{C\left(v_j,\frac{r}{2}\right)}} \tilde{z}>u\right\}}=U(u).
\end{align*}
\end{proof}

	\section{FWER control and power consistency }{\label{sec: proof FWER}}
	\begin{proof}[Proof of Theorem~\ref{theorem:FWER}]
		We start with the power. Since our algorithm erases cubes of side length $r=2B+\delta$ (step $2$ in \Cref{sec:algorithm}) and $\delta\leq \frac{r}{2}$ (because $\delta\leq 2B$, see \Cref{appendix: delta zero estimation}), there is at most one element from $T$ in each $\delta$-neighborhood of the centers.  Then for all $u\in \mathbb{R}$,  
			\begin{equation*}
			\text{Power}(u)=\mathbb{E}\left[\frac{1}{N}\sum_{i=1}^{N}\mathbbm{1}_{\{T(u)\cap C(\tau_i,\delta)\neq\emptyset\}}\right]=\mathbb{E}\left[\frac{1}{N}\#\left\{T(u)\cap\mathbb{B}_{1}^{\delta}\right\}\right].
			\end{equation*}
		We shall prove that
		\begin{equation*}
			\text{Power}(u_{\rm{Bon}})=\mathbb{E}\left[\frac{1}{N}\#\left\{T\left(u_{\rm{Bon}}\right)\cap\mathbb{B}_{1}^{\delta}\right\}\right]\underset{L\rightarrow\infty}{\longrightarrow}1.
		\end{equation*}
		On the event $ E_L =   E^{(1)}_L\cap E^{\text{Bon}}_L $, it holds that $\#\left\{T\left(u_{\rm{Bon}}\right)\cap\mathbb{B}_{1}^{\delta}\right\}=N$.
		By Claims \ref{claim: the first N points in the algorithm are objects} and~\ref{claim: BH algorithm have more then N steps}, $$\mathbb{P}\left[\left(E_L\right)^{\complement}\right]=\mathbb{P}\left[\left(E_L^{(1)}\right)^{\complement}\cup\left(E_L^{\text{Bon}}\right)^{\complement}\right]\lesssim 2^N\cdot L^de^{-a_L^2} + N\cdot\frac{1}{L^{4d}}\xrightarrow{L\rightarrow\infty}0.$$
		Hence, we get
		\begin{align*}
			\text{Power}(u_{\rm{Bon}})=\mathbb{E}\left[\frac{1}{N}\#\left\{T\left(u_{\rm{Bon}}\right)\cap\mathbb{B}_{1}^{\delta}\right\}\right]&\geq \mathbb{E}\left[\frac{1}{N}\#\left\{T\left(u_{\rm{Bon}}\right)\cap\mathbb{B}_{1}^{\delta}\right\}\cdot\mathbbm{1}_{E_L}\right]&\\= \mathbb{E}\left[\mathbbm{1}_{E_L}\right]=\mathbb{P}\left[E_L\right]\xrightarrow{L\rightarrow\infty}1.
		\end{align*} 
	 Recall that for all $u\in \mathbb{R}$,  $U(u)=\sum_{j=1}^{M_L}\mathbbm{1}_{\left\{\max_{\overline{C\left(v_j,\frac{r}{2}\right)}} \tilde{z}>u\right\}}$ and that the test function \ref{eq: test value} is  defined by $p(u) = \mathbb{P}\left[ \max_{\overline{C\left(0,\frac{r}{2}\right)}}\tilde{z}>u\right]$. Note that due to stationarity of $\tilde{z}$ we have that $$p(u) = \frac{\mathbb{E}\left[U(u)\right]}{M_L} .$$
	Next we show that the FWER is controlled i.e. 
	\begin{equation*}
		\limsup_{L\rightarrow\infty}\text{FWER}(u_\text{Bon})=\limsup_{L\rightarrow\infty}\mathbb{P}[V(u_{\rm{Bon}})\geq 1]\leq \alpha.
	\end{equation*}  
	 Now, applying Markov's inequality and using Claim \ref{claim: V less U} and the trivial bound $V(u)\leq M_L$, we have that, 
	\begin{align*}
		\mathbb{P}[V(u_{\rm{Bon}})\geq 1]\leq \mathbb{E}[V(u_{\rm{Bon}})]&=  \mathbb{E}[V(u_{\rm{Bon}})\cdot\mathbbm{1}_{E_L^{(1)}}]+\mathbb{E}[V(u_{\rm{Bon}})\cdot\mathbbm{1}_{(E_L^{(1)})^{\complement}}]\\ &\leq \mathbb{E}[U(u_{\rm{Bon}})\cdot\mathbbm{1}_{E_L^{(1)}}]+M_L\cdot\mathbb{E}[\mathbbm{1}_{(E_L^{(1)})^{\complement}}]\\&\leq \mathbb{E}[U(u_{\rm{Bon}})]+M_L\cdot\mathbb{E}[\mathbbm{1}_{(E_L^{(1)})^{\complement}}]\\&\leq p(u_{\rm{Bon}})\cdot M_L+ c M_L 2^N L^de^{-a_L^2}\\&=\alpha +c M_L 2^N L^de^{-a_L^2},
	\end{align*} 
	where the  fifth inequality follows from \Cref{claim: the first N points in the algorithm are objects}.. 
	Taking $\limsup_{L\rightarrow\infty}$ concludes the proof by Assumptions  \ref{assumption: limit of the number of objects} and \ref{assumption: limit of snr}.
	\end{proof}
	\section{FDR control and power consistency}\label{sec: proof FDR}

	\begin{proof}[Proof of Theorem~\ref{theorem: fdr}]
		The proof of power consistency is the same as in Appendix~\ref{sec: proof FWER}.  For the FDR, we shall use subsequences of $L$, so we denote $u_{\rm{BH}}=u_{L}^{\rm{BH}}$, $V(u)=V_L(u)$, $W(u)=W_L(u)$, and $A'=\frac{A}{B^d}$. We need to prove that 
	\begin{equation*}
		\limsup_{L\rightarrow\infty}\text{FDR}(u_{L}^{\rm{BH}})=\limsup_{L\rightarrow\infty}\mathbb{E}\left[\frac{V_L(u_{L}^{\rm{BH}})}{V_L(u_{L}^{\rm{BH}})+W_L(u_{L}^{\rm{BH}})}\right]\leq \alpha.
	\end{equation*}
	Without the loss of generality, assume $L\in\mathbb{N}$. 
	Set $E_L=E_L^{(1)}\cap E_L^{\text{BH}}$ to be the event where the first $N$ points in~$T$ are in~$\mathbb{B}_1^{\delta}$, and that the BH algorithm (defined in \Cref{sec:theorem}) stops after more than $N$ steps i.e. $k>N$. By Claims \ref{claim: the first N points in the algorithm are objects}, \ref{claim: BH algorithm have more then N steps} and Assumption~\ref{assumption: limit of the number of objects} in Section~\ref{sec:theorem}, there exists $L_0$ such that for all $L>L_0$, 
	\begin{equation}\label{eq: L_0 def}
		\mathbb{P}[{E_{L}}^{\complement}]\lesssim 2^NL^de^{-a_L^2}+N\cdot\frac{1}{L^{4d}}\lesssim \frac{1}{L^2}, \quad \frac{N}{M_L}\geq \frac{A'(\frac{r}{2})^d}{2}.
	\end{equation}
	Now, consider the event $$E_{L_0}=\bigcap_{L\geq L_0}E_{L},$$ and note that
	\begin{align*}
		\mathbb{P}[{E_{L_0}}^{\complement}]\leq \sum_{n\geq L_0}\mathbb{P}[{E_{L}}^{\complement}]\lesssim\sum_{L\geq L_0}\frac{1}{L^2}\lesssim \frac{1}{L_0-1}\xrightarrow{L_0\rightarrow\infty}0.
	\end{align*} 
	Throughout the proof, we will rely on the ergodicity of the process $z$ (defined in \eqref{eq: the model}) which implies that for all $u\in\mathbb{R}$, 
	\begin{equation}\label{eq: ergodicty}
		\frac{U_L(u)}{M_L}=\frac{\sum_{j=1}^{M_L}\mathbbm{1}_{I_j(u)}}{M_L}\xrightarrow{L\rightarrow\infty}p(u) \text{ a.s}.
	\end{equation}
	 Set  \begin{equation}\label{def: G}
	 	G=\bigcap_{q\in\mathbb{Q}}G_q,
	 \end{equation}
 where $G_q=\{\omega\in \Omega \mid \frac{(U_L(q))(\omega)}{M_L}\xrightarrow{L\to\infty}p(q)\}$ and $\Omega$ is our probability space for the process $z$.  Then, by \eqref{eq: ergodicty}, for all $q\in\mathbb{Q}$, $\mathbb{P}[G_q]=1$,
	which implies $\mathbb{P}\left[G\right]=1$.  
	Set $\tilde{E}_{L_0}=E_{L_0}\cap G$ and note that $\mathbb{P}[\tilde{E}_{L_0}]=\mathbb{P}[E_{L_0}]$. We have,
		\begin{align*}
		\text{FDR}(u_{L}^{\rm{BH}})&=\mathbb{E}\left[\frac{V_L(u_{L}^{\rm{BH}})}{V_L(u_{L}^{\rm{BH}})+W_L(u_{L}^{\rm{BH}})}\cdot \mathbbm{1}_{\tilde{E}_{L_0}}\right]+\mathbb{E}\left[\frac{V_L(u_{L}^{\rm{BH}})}{V_L(u_{L}^{\rm{BH}})+W_L(u_{L}^{\rm{BH}})}\cdot \mathbbm{1}_{{\tilde{E}_{L_0}}^{\complement}}\right]\\ &\leq \mathbb{E}\left[\frac{V_L(u_{L}^{\rm{BH}})}{V_L(u_{L}^{\rm{BH}})+W_L(u_{L}^{\rm{BH}})}\cdot \mathbbm{1}_{\tilde{E}_{L_0}}\right]+\mathbb{P}\left[{\tilde{E}_{L_0}}^{\complement}\right].
	\end{align*}
	Now, for all $L\geq L_0$
	\begin{equation*}
		\frac{V_L(u_{L}^{\rm{BH}})}{V_L(u_{L}^{\rm{BH}})+W_L(u_{L}^{\rm{BH}})}\cdot \mathbbm{1}_{\tilde{E}_{L_0}}=\frac{V_L(u_{L}^{\rm{BH}})\cdot\mathbbm{1}_{\tilde{E}_{L_0}}}{V_L(u_{L}^{\rm{BH}})\cdot\mathbbm{1}_{\tilde{E}_{L_0}}+N}.
	\end{equation*}
	Hence
	\begin{align*}
		&\mathbb{E}\left[\frac{V_L(u_{L}^{\rm{BH}})}{V_L(u_{L}^{\rm{BH}})+W_L(u_{L}^{\rm{BH}})}\cdot \mathbbm{1}_{\tilde{E}_{L_0}}\right]= \mathbb{E}\left[\frac{V_L(u_{L}^{\rm{BH}})\cdot\mathbbm{1}_{\tilde{E}_{L_0}}}{V_L(u_{L}^{\rm{BH}})\cdot\mathbbm{1}_{\tilde{E}_{L_0}}+N}\right]\\&\leq \frac{\mathbb{E}[V_L(u_{L}^{\rm{BH}})\cdot\mathbbm{1}_{\tilde{E}_{L_0}}]}{\mathbb{E}[V_L(u_{L}^{\rm{BH}})\cdot\mathbbm{1}_{\tilde{E}_{L_0}}]+N},
	\end{align*}
	where the last inequality follows from Jensen's inequality. By claim \ref{claim: V less U}, we have that for all $L\geq L_0$, on the event $\tilde{E}_{L_0}$, $V_L(u_L^{\text{BH}})\leq U_L(u_L^{\text{BH}})$. By monotonicity of the map $x\mapsto \frac{x}{x+a}\text{ for }x,a>0$,
	$$\text{FDR}(u_{L}^{\rm{BH}})\leq \frac{\mathbb{E}[U_L(u_{L}^{\rm{BH}})\cdot\mathbbm{1}_{\tilde{E}_{L_0}}]}{\mathbb{E}[U_L(u_{L}^{\rm{BH}})\cdot\mathbbm{1}_{\tilde{E}_{L_0}}]+N}+\mathbb{P}\left[{\tilde{E}_{L_0}}^{\complement}\right].$$
	%
	The latter expression is connected to the test value \eqref{eq: test value} as for a  deterministic threshold $u$, $\mathbb{E}[U_L(u)]/ M_L=p(u)$. Though $u^{\text{BH}}_L$ is random, the next lemma overcomes the randomness by bounding it from below with a random sequence that converges to a deterministic threshold.  
	\begin{lemma}\label{lemma: uL to u star}
		Let $u^*\in\mathbb{R}^{+}$ be the smallest solution of $p(u)=\frac{\alpha A'(\frac{r}{2})^d}{1-\alpha}$.
		Then, there exists a random sequence $u_L$ such that for all $L\geq L_0$ and $\omega\in \tilde{E}_{L_0}$
		\begin{itemize}
			\item $u_L(\omega)\leq u_L^{\text{BH}}(\omega)$
			\item 	$u_L(\omega)\xrightarrow{L\rightarrow\infty}u^*$.
		\end{itemize} 
	\end{lemma} 
	\begin{proof}[The Proof of Lemma \ref{lemma: uL to u star} is in \Cref{proof: of lemma uL}]
\end{proof}
	This lemma, combined with the monotonicity of $U_L$, implies that 
	$$\frac{\mathbb{E}[U_L(u_{L}^{\rm{BH}})\cdot\mathbbm{1}_{\tilde{E}_{L_0}}]}{\mathbb{E}[U_L(u_{L}^{\rm{BH}})\cdot\mathbbm{1}_{\tilde{E}_{L_0}}]+N}\leq 	\frac{\mathbb{E}[U_L(u_{L})\cdot\mathbbm{1}_{\tilde{E}_{L_0}}]}{\mathbb{E}[U_L(u_{L})\cdot\mathbbm{1}_{\tilde{E}_{L_0}}]+N}.$$
	By Egorov's Theorem \cite{kolmogorov1975introductory}, there exists $D_{L_0}\subset \tilde{E}_{L_0}$ such that $\mathbb{P}[\tilde{E}_{L_0}\backslash D_{L_0}]<\frac{1}{L_0}$ and $$u_L\cdot\mathbbm{1}_{D_{L_0}}\xrightarrow{uni}u^*\cdot\mathbbm{1}_{D_{L_0}},$$
	which yields that there exists $L_1>L_0$ such that for all $L\geq L_1$ and $\omega\in D_{L_0}$
	$$u_L(\omega)>u^*-\frac{1}{L_0}.$$
	It follows that
	\begin{align*}
		\frac{\mathbb{E}[U_L(u_{L})\cdot\mathbbm{1}_{\tilde{E}_{L_0}}]}{\mathbb{E}[U_L(u_{L})\cdot\mathbbm{1}_{\tilde{E}_{L_0}}]+N}&=\frac{\mathbb{E}[U_L(u_{L})\cdot\mathbbm{1}_{D_{L_0}}]+\mathbb{E}[U_L(u_{L})\cdot\mathbbm{1}_{\tilde{E}_{L_0}\backslash D_{L_0}}]}{\mathbb{E}[U_L(u_{L})\cdot\mathbbm{1}_{D_{L_0}}]+\mathbb{E}[U_L(u_{L})\cdot\mathbbm{1}_{\tilde{E}_{L_0}\backslash D_{L_0}}]+N}\\&\leq \frac{\mathbb{E}[U_L(u_{L})\cdot\mathbbm{1}_{D_{L_0}}]+\mathbb{E}[U_L(u_{L})\cdot\mathbbm{1}_{\tilde{E}_{L_0}\backslash D_{L_0}}]}{\mathbb{E}[U_L(u_{L})\cdot\mathbbm{1}_{D_{L_0}}]+N}\\&= 
		\frac{\mathbb{E}[U_L(u_{L})\cdot\mathbbm{1}_{D_{L_0}}]}{\mathbb{E}[U_L(u_{L})\cdot\mathbbm{1}_{D_{L_0}}]+N}+\frac{\mathbb{E}[U_L(u_{L})\cdot\mathbbm{1}_{\tilde{E}_{L_0}\backslash D_{L_0}}]}{\mathbb{E}[U_L(u_{L})\cdot\mathbbm{1}_{D_{L_0}}]+N}
		\\&\overset{\star}{\leq}  
		\frac{\mathbb{E}[U_L(u^*-\frac{1}{L_0})\cdot\mathbbm{1}_{D_{L_0}}]}{\mathbb{E}[U_L(u^*-\frac{1}{L_0})\cdot\mathbbm{1}_{D_{L_0}}]+N}+\frac{M_L}{N}\cdot\mathbb{P}[\tilde{E}_{L_0}\backslash D_{L_0}]
		\\&
		\leq 	 \frac{\mathbb{E}[U_L(u^*-\frac{1}{L_0})]/M_L}{\mathbb{E}[U_L(u^*-\frac{1}{L_0})]/M_L+N/M_L}+\frac{M_L}{N}\cdot\mathbb{P}[\tilde{E}_{L_0}\backslash D_{L_0}]\\&= \frac{p(u^*-\frac{1}{L_0})}{p(u^*-\frac{1}{L_0})+N/M_L}+\frac{M_L}{N}\cdot\mathbb{P}[\tilde{E}_{L_0}\backslash D_{L_0}],
	\end{align*}
	where in $\star$ we trivially bound $U_L\leq M_L$. 
	To conclude,  for all $L\geq L_1\geq L_0$
	\begin{equation*}
		\text{FDR}(u_{L}^{\rm{BH}})\leq  \frac{p(u^*-\frac{1}{L_0})}{p(u^*-\frac{1}{L_0})+N/M_L}+\frac{M_L}{N}\cdot\mathbb{P}[\tilde{E}_{L_0}\backslash D_{L_0}]+ \mathbb{P}\left[\left(\tilde{E}_{L_0}\right)^{\complement}\right].
	\end{equation*}
	Hence, 
	\begin{equation*}
		\limsup_{L\to\infty}\text{FDR}(u_{L}^{\rm{BH}})\leq  \frac{p(u^*-\frac{1}{L_0})}{p(u^*-\frac{1}{L_0})+A'(\frac{r}{2})^d}+\frac{1}{A'(\frac{r}{2})^d}\cdot\frac{1}{L_0}+ \frac{c}{L_0},
	\end{equation*}
        where $c$ is a positive constant. 
	Taking $L_0\to\infty$, by continuity of $p$  (see \Cref{claim: contiuity of p}), we get
	\begin{equation*}
		\limsup_{L\to\infty}\text{FDR}(u_{L}^{\rm{BH}})\leq  \frac{p(u^*)}{p(u^*)+A'(\frac{r}{2})^d}=\frac{\frac{\alpha A'(\frac{r}{2})^d}{1-\alpha}}{\frac{\alpha A'(\frac{r}{2})^d}{1-\alpha}+A'(\frac{r}{2})^d}=\alpha.
	\end{equation*}
\end{proof}
 \setcounter{subsection}{1}\subsection {Proof of Lemma \ref{lemma: uL to u star} }\label{proof: of lemma uL}
	
		First, we will prove that on $\tilde{E}_{L_0}$,  every convergent subsequence of $u_L^{\text{BH}}$  converge to a solution of the inequality $p(v)\leq\frac{A'(\frac{r}{2})^d\alpha}{1-\alpha}$. 
		Recall that on $\tilde{E}_{L_0}$, The first $N$ points in the candidate set $T$  are in  $\mathbb{B}^1_\delta$, and the Benyamini-Hochberg procedure (defined in \cref{sec:theorem}) stops after more then $N$ steps. As a result, the number of detected objects, denoted by $k$, equal to the number of objects plus the number of false positives. The latter, combined with the definition of the BH-threshold $u_{L}^{\rm{BH}}$ (see \Cref{sec:theorem}), gives for all $L\geq L_0$, \begin{equation*}
			p(u_{L}^{\rm{BH}})=\frac{k\alpha}{M_L}=\frac{\left(N+V_L\left(u_{L}^{\rm{BH}}\right)\right)\alpha}{M_L}.
		\end{equation*}
		Let $\omega\in \tilde{E}_{L_0}$, and note that $u_L^{\text{BH}}(\omega)$ is bounded. Indeed, by \eqref{eq: L_0 def} for all $L\geq L_0$, $\frac{N}{M_L}>\frac{A'(\frac{r}{2}))^d}{2}$ and so, 
		\begin{equation}
			p\left(u_L^{\text{BH}}(\omega)\right)=\frac{\left(N+\left(V_{L}((\omega))\right)\left(\omega\right)\right)\alpha}{M_{L}}> \frac{A'\left(\frac{r}{2}\right)^d\alpha}{2}.
		\end{equation}
		 Moreover, by definition, $k\leq m_L < M_L$  and hence
		\begin{equation}
			p\left(u_L^{\text{BH}}(\omega)\right)=\frac{k\alpha}{M_L}<\frac{M_L\alpha}{M_L}\leq1.
		\end{equation} 
		Since $p$ is monotone decreasing, $u_L^{\text{BH}}(\omega)$ is bounded.
		Consider $u_{L_j}^{\text{BH}}(\omega) \xrightarrow{j\rightarrow\infty}v\in\mathbb{R}^{+}$, a convergent subsequence of $u_L^{\text{BH}}(\omega)$. Then, for every $\delta>0$ there exists a large enough $j_0$ such that for every $j>j_0$, 
		\begin{equation*}
			v-\delta<u_{L_j}^{\text{BH}}(\omega).
		\end{equation*}
		Let $\delta>0$, and choose $\delta_1 < \delta$ such that $v - \delta_1 \in \mathbb{Q}$. 
		By \Cref{claim: V less U} and the monotonicity of $U_L$, we get for all $j\geq j_0$,
		\begin{align}\label{eq: inequality}
			p(u_{L_j}^{\text{BH}}(\omega))&=\frac{\left(N+\left(V_{L_j}(u_{L_j}^{\text{BH}}(\omega))\right)\left(\omega\right)\right)\alpha}{M_{L_j}}\\ \nonumber&\leq \frac{\left(N+\left(U_{L_j}(u_{L_j}^{\text{BH}}(\omega))\right)\left(\omega\right)\right)\alpha}{M_{L_j}}\\ \nonumber &\leq \frac{\left(N+\left(U_{L_j}\left(v-\delta_1\right)\right)\left(\omega\right)\right)\alpha}{M_{L_j}}.
		\end{align}
	Since $\omega\in G$ (defined in \eqref{def: G}), taking the limit $j\rightarrow \infty$ in~\eqref{eq: inequality} and using the continuity of $p$ yields 
		\begin{equation*}
			p(v)\leq \left(\left(\frac{r}{2}\right)^dA'+p(v-\delta_1)\right)\alpha.
		\end{equation*}
	 Using again the monotonicity of $p$, we get that 	
		$$p(v)\leq \left(\left(\frac{r}{2}\right)^dA'+p(v-\delta)\right)\alpha.$$ 
	 Taking limit over $\delta\rightarrow{0}$ in the last inequality implies that $v$ satisfies $$p(v)\leq	\left(\left(\frac{r}{2}\right)^dA'+p(v)\right)\alpha,$$ 
        hence
        $$p(v)\leq\frac{\alpha A'\left(\frac{r}{2}\right)^d}{1-\alpha}.$$
		Set $u_L=\min\{u_L^{\text{BH}},u^*\}$, where~$u^*$ is the smallest solution for $p(u)=\frac{\alpha A'\left(\frac{r}{2}\right)^d}{1-\alpha}$ and
		let $\omega\in \tilde{E}_{L_0}$. We want to prove that $u_L(\omega)\xrightarrow{L\rightarrow\infty}u^*$.
		Let $u_{L_l}(\omega)\xrightarrow{j\rightarrow\infty} s\in\mathbb{R}^+$. Since $ u_{L_l}^{\text{BH}}(\omega) $ is bounded, there exists a convergent  subsequence $u_{L_{l_i}}^{\text{BH}}(\omega)\xrightarrow{i\rightarrow\infty}v$ where $p(v)\leq\frac{A'(\frac{r}{2})^d\alpha}{1-\alpha}$. Now, 
		\begin{align*}
			&u_{L_{l_i}}(\omega)=\min\left\{u^{BH}_{L_{l_i}}(\omega),u^*\right\}\xrightarrow{i\rightarrow\infty}\min\{v,u^*\}=u^*,\\\nonumber
			\end{align*} 
        which implies $s=u^*$. We proved that every convergent subsequence of $ u_{L}(\omega) $ converges to the same deterministic $u^*$, and therefore $ u_{L}(\omega) \xrightarrow{L\to\infty} u^*$.\\
		To conclude, we found a random sequence $u_L$ such that for all $L\geq L_0$ and $\omega\in \tilde{E}_{L_0}$, $u_L(\omega)\leq u_L^{\text{BH}}(\omega)$ and $u_L(\omega)\xrightarrow{L\rightarrow\infty}u^*$ as desired. 
  \section{Estimating \texorpdfstring{$\delta$}{delta} which satisfies Definition~\ref{def: localization property}}{} \label{appendix: delta zero estimation}
  Recall \Cref{def: localization property}:
		we say that $\delta > 0$ satisfies the localization property if for all $1\leq i\leq N$, there exists some representative $q_i\in C(\tau_i,\delta)$ such that $\forall t\in \overline{C(0,L)}\backslash \left(\cup_{j=1}^{i-1}C(\tau_j,2B) \cup_{j=i}^{N}C(\tau_j,\delta)\right)$
		\begin{equation*}
			S^x(q_i)-S^x(t)\geq \rho ||a_i||^2,
		\end{equation*}
		where $	S^x = \sum_{j=1}^M \left(x\ast\tilde{\psi}_j\right)^2$, 
		 and $\rho$ is a positive constant.
     Clearly, such $\delta$ depends on the given  $\{\psi_j\}_{j=1}^{M}$, the unknown coefficients $\{a_i\}_{i=1}^{N}$ and locations of the objects, all are defined in \Cref{sec:setup}. Nonetheless, in all cases,  $\delta=2B$ satisfies \Cref{def: localization property}, hence such $\delta$s always exists.  Indeed, for $\delta=2B$, it holds that $S^x(t)=0$ for all $t\in\mathbb{B}_0^\delta$, and by  \Cref{claim: centers norms}, it holds that $S^x(\tau_i)=||a_i||^2$ for all $1\leq i\leq N$. Therefore, taking $q_i=\tau_i$ suffices. Our goal is to find $\delta$ satisfying Definition~\ref{def: localization property} which is smaller than~$2B$. To that end, we will  develop an inequality whose solutions satisfy \Cref{def: localization property}, see \Cref{claim: inequality yields delta property}, and use it throughout \Cref{sec: numerical results} in order to validate that a given value of~$\delta$ can be used.

     First we note that $S^x$ can be written as a bi-linear form involving the inner product of the translated basis functions. Set the matrix $D(t,\tau_i,\tau_l)\in\mathbb{R}^{M\times M}$ to be, for all $1\leq k,s \leq M$, 
     $$ D_{k,s}(t,\tau_i,\tau_l) = \sum_{j=1}^M\int_{C\left(0,B\right)}\psi_k(u)\psi_j\left(u-\left(t-\tau_i\right)\right)du\int_{C\left(0,B\right)}\psi_s(u')\psi_j\left(u'-\left(t-\tau_l\right)\right)du'.$$   $ D(t,\tau_i,\tau_l) $ is a  real symmetric matrix, which depends on $(t-\tau_i,t-\tau_l)$. 
     \begin{claim}
            \begin{equation}\label{eq: S^x as bi linear form}
     	S^x(t) =\sum_{i,l=1}^N a^T_i D(t,\tau_i,\tau_l)a_l.
     \end{equation} 
     \end{claim}
 	\begin{proof}
 	Denote  $ \psi^i_j(t)=\psi_j(t-\tau_i) $. Then,
 		\begin{align*}
 			S^x(t) &= \sum_{j=1}^M \left(x\ast\tilde{\psi_j}\left(t\right)\right)^2 =\sum_{j=1}^M \left(\sum_{i=1}^Nx_i\ast\tilde{\psi_j}\left(t\right)\right)^2\\\nonumber
 			&=\sum_{j=1}^M \left(\sum_{i=1}^N\sum_{k=1}^Ma_{ik}\psi^i_k \ast\tilde{\psi_j}\left(t\right)\right)^2\\\nonumber
 			&=\sum_{j=1}^M \left(\sum_{i=1}^N\sum_{k=1}^Ma_{ik}\int_{\mathbb{R}^d}\psi^i_k(y)\tilde{\psi}_j\left(t-y\right)dy\right)^2\\\nonumber
 			&=\sum_{j=1}^M \left(\sum_{i=1}^N\sum_{k=1}^Ma_{ik}\int_{\mathbb{R}^d}\psi_k(y-\tau_i)\psi_j\left(y-t\right)dy\right)^2\\\nonumber
 			&=_{\left[y-\tau_i=u\right]}\sum_{j=1}^M \left(\sum_{i=1}^N\sum_{k=1}^Ma_{ik}\int_{C\left(0,B\right)}\psi_k(u)\psi_j\left(u-\left(t-\tau_i\right)\right)du\right)^2\\\nonumber
 			&=\sum_{j=1}^M \sum_{i,l=1}^N\sum_{k,s=1}^Ma_{ik}a_{ls}\int_{C\left(0,B\right)}\psi_k(u)\psi_j\left(u-\left(t-\tau_i\right)\right)du\int_{C\left(0,B\right)}\psi_s(u')\psi_j\left(u'-\left(t-\tau_l\right)\right)du'.\\\nonumber
 			&=\sum_{i,l=1}^N\sum_{k,s=1}^Ma_{ik}a_{ls}\sum_{j=1}^M\int_{C\left(0,B\right)}\psi_k(u)\psi_j\left(u-\left(t-\tau_i\right)\right)du\int_{C\left(0,B\right)}\psi_s(u')\psi_j\left(u'-\left(t-\tau_l\right)\right)du'.
 			\\&=\sum_{i,l=1}^N \sum_{k,s=1}^N a_{ik}D_{k,s}(t,\tau_i,\tau_l)a_{ls}  =\sum_{i,l=1}^N a^T_iD(t,\tau_i,\tau_l)a_l.
 		\end{align*}
 	\end{proof}
		Denote by
		\begin{align*}
			\lambda_{\min}(t,\tau_i,\tau_i)&=\lambda_{\min} \left(D(t,\tau_i,\tau_i)\right),\\\nonumber
			\lambda_{\max}(t,\tau_i,\tau_l)&= \max\left\{\lvert\lambda_{\min}\left(D(t,\tau_i,\tau_l)\right)\rvert,\lvert\lambda_{\max}\left(D(t,\tau_i,\tau_l)\right)\rvert\right\},
		\end{align*}
where $\lambda_{\min},\lambda_{\max}$ are the minimal and maximal eigenvalue respectively.
  
   \begin{claim}\label{claim: bounds on Sx}
     Let $\delta\in [0,2B]$, and assume $||\tau_j-\tau_l||_\infty\geq B+\frac{3}{2}\delta$ for all $j\neq l$. Denote $I(t)=\{1\leq j\leq N \mid ||\tau_j-t||_\infty< B\}$. Then, 
     \begin{enumerate}
         \item $\forall 1\leq i\leq N$, $\forall t\in C(\tau_i,\delta),\;\; S^x(t)\geq ||a_i||^2\lambda_{\min}(t,\tau_i,\tau_i).$
         \item $\forall t'\in\mathbb{B}_0^{\delta},\;\; S^x(t')\leq \sum_{j,l\in I(t')}||a_j||\cdot||a_l||\lambda_{\max}(t',\tau_j,\tau_l)$.
     \end{enumerate}
 \end{claim}
  \begin{proof}
  For the first part, let $t\in C(\tau_i,\delta)$, then
    \begin{align}
		S^x(t) & =\sum_{j,l=1}^N a^T_jD(t,\tau_j,\tau_l)a_l \\\nonumber
		& =a_i^{T}D(t,\tau_i,\tau_i)a_i\geq ||a_i||^2\lambda_{\min}(t,\tau_i,\tau_i).
  \end{align}
  The second equality in the latter expression is due to the fact that if $t\in C(\tau_i,\delta)$ and the centers are  more than $B+\frac{3\delta}{2}$ apart, then all the other terms in the sum of $S^x$ are zero.  
    The last inequality in the latter expression is a known bound in linear algebra \cite{hoffman1971linear}.
     
     For the second part, let $t'\in\mathbb{B}_0^{\delta}$. Then,
    \begin{align}
		S^x(t') & =\sum_{j,l\in I(t')} a^T_jD(t',\tau_j,\tau_l)a_l \\\nonumber
		& \leq \sum_{j,l\in I(t')}\norm{a_j}\norm{a_l}\lambda_{\max}(t',\tau_j,\tau_l),
  \end{align}
 where  the last inequality is again a known bound from linear algebra \cite{hoffman1971linear}.
  \end{proof}
  By \Cref{claim: bounds on Sx}, we have  $\forall t\in C(\tau_i,\delta)$ and  $\forall t'\in \mathbb{B}_{0}^{\delta}$,
  $$S^x(t)-S^x(t')\geq ||a_i||^2\lambda_{\text{min}}(t,\tau_i,\tau_i)-\sum_{j,l\in I(t')}\norm{a_j}\norm{a_l}\lambda_{\max}(t',\tau_j,\tau_l).$$
     As mentioned above, $D$ is a function of $(t-\tau_j,t-\tau_l)$, which makes $D$ 
    translation invariant with respect to the object centers and $t$. That is, for all $v\in\mathbb{R}^d$ we have $D\left(t,\tau_j,\tau_l\right) = D(t-v,\tau_j-v,\tau_l-v)$, which implies the same for the eigenvalues of~$D$. 
    In particular, for $t\in C(\tau_i,\delta)$, $$\lambda_{\text{min}}(t,\tau_i,\tau_i)=\lambda_{\text{min}}(t-\tau_i,0,0),$$ and for $t'\in \mathbb{B}_0^{\delta}$ $$\lambda_{\max}(t',\tau_j,\tau_l)=\lambda_{\max}(0,\tau_j-t',\tau_l-t').$$ 
     To conclude, $$S^x(t)-S^x(t')\geq \norm{a_i}^2\lambda_{\min}(t-\tau_i,0,0) - \sum_{j,l\in I(t')}\norm{a_j}\norm{a_l}\lambda_{\max}(0,\tau_j-t',\tau_l-t'), $$ where $I(t)$ is defined on \Cref{claim: bounds on Sx}. Now, taking the supremum with respect to $t,t'$ and all possible locations of objects yields  the following definition.
    \begin{definition}\label{def: g(delta)}
    \begin{align*}
       g(\delta) = \sup_{t\in C(0,\delta)}\lambda_{\min}(t,0,0)-\sup_{\substack{\norm{\tau_l-\tau_j}\geq B+\frac{3}{2}\delta\\\tau_l\notin C(0,\delta)}}\;\sum_{l,j\in I(0)}\lambda_{\max}(0,\tau_l,\tau_j)
    \end{align*}
    for $0\leq \delta\leq 2B.$
    \end{definition}

    \begin{claim}\label{claim: inequality yields delta property}
       Assume $\forall 1\leq j\neq l\leq N$ that  $||\tau_j-\tau_l||\geq B+\frac{3\delta}{2}$. Then, given $\delta > 0$ such that $g(\delta) > 0$, implies that $\delta$  satisfies \Cref{def: localization property}. 
    \end{claim}
    \begin{proof}
        Choose $\epsilon>0$ such that $g(\delta)-\epsilon>0$. By \Cref{def: g(delta)}, $\exists q_1\in C(0,\delta)$ for which $\lambda_{\text{min}}(q_1,0,0)\geq \sup_{t\in C(0,\delta)}\lambda_{\min}(t,0,0)-\epsilon$. Using the shift invariant property of $D$ and \Cref{def: g(delta)}:
        \begin{align}
        	\forall t'\in\mathbb{B}_0^{\delta}, \left(\lambda_{\text{min}}(q_1,\tau_1,\tau_1)-\sum_{j,l\in I(t')}\lambda_{\max}(t',\tau_j,\tau_l)\right) \geq g(\delta)-\epsilon
        \end{align} 
    Assume without loss of generality that $||a_1||\geq\dots\geq||a_N||$. Then,
        by \Cref{claim: bounds on Sx}
        	 \begin{align*}
            S^x(q_1)-S^x(t')&\geq||a_1||^2\lambda_{\text{min}}(q_1,\tau_1,\tau_1)-\sum_{j,l\in I(t')}\norm{a_j}\norm{a_l}\lambda_{\max}(t',\tau_j,\tau_l)\\ & \geq  ||a_1||^2\lambda_{\text{min}}(q_1,\tau_1,\tau_1)-||a_1||^2\sum_{j,l\in I(t')}\lambda_{\max}(t',\tau_j,\tau_l)\\ &= ||a_1||^2\left(\lambda_{\text{min}}(q_1,\tau_1,\tau_1)-\sum_{j,l\in I(t')}\lambda_{\max}(t',\tau_j,\tau_l)\right)\\& \geq ||a_1||^2(g(\delta)-\epsilon).
        \end{align*}
        Next, if $t'\in \overline{C(0,L)}\backslash  C(\tau_1,2B) \cup _{j=2}^N C(\tau_j,\delta)$, then $||\tau_1-t'||\geq B$, which implies $1\notin I(t').$ where $I(t')$ defined in \Cref{claim: bounds on Sx}. By \Cref{claim: bounds on Sx}, the shift invariant property of $D$, and \Cref{def: g(delta)}, $\exists q_2\in C(\tau_2,\delta)$ such that $\forall t'\in \overline{C(0,L)}\backslash  C(\tau_1,2B) \cup _{j=2}^N C(\tau_j,\delta)$,
        \begin{align*}
            S^x(q_2)-S^x(t')&\geq||a_2||^2\lambda_{\text{min}}(q_2,\tau_2,\tau_2)-\sum_{j,l\in I(t')}\norm{a_j}\norm{a_l}\lambda_{\max}(t',\tau_j,\tau_l)\\ & \geq  ||a_2||^2\lambda_{\text{min}}(q_2,\tau_2,\tau_2)-||a_2||^2\sum_{j,l\in I(t')}\lambda_{\max}(t',\tau_j,\tau_l)\\ &= ||a_2||^2\left(\lambda_{\text{min}}(q_2,\tau_2,\tau_2)-\sum_{j,l\in I(t')}\lambda_{\max}(t',\tau_j,\tau_l)\right)\\& \geq ||a_2||^2(g(\delta)-\epsilon).
        \end{align*}
        We can repeat this argument for $3\leq i \leq N$, and get the desired property, which completes the proof. 
    \end{proof}

	Recall that our algorithm requires as input a $\delta$ that satisfies  \Cref{def: localization property} (to ensure the statistical guarantees outlined in \Cref{theorem:FWER,theorem: fdr}). When processing numerical data, as discussed in \Cref{sec: numerical results}, the basis functions spanning the objects are represented on a grid that scales with the size of the objects (denoted by $B$). Following the assertion made in \Cref{claim: inequality yields delta property}, we seek to determine such a value for $\delta$, by numerically solving the inequality $g(\delta)>0$. To that end, we iterate over $\delta$, starting  with the smallest possible value which is computationally feasible, for example, the minimum distance between two points on the grid. In each iteration, in order to estimate $g(\delta)$, we  computes $\lambda_{\max},\lambda_{\min}$ for all possible centers $\{\tau_i\}_i$ and points $t$  on the grid.  This iterative process persists until $g(\delta)>0$ invariably concluding at  $\delta=2B$ since $g(2B)>0$ (as explained earlier in this section).

\end{appendices}

\bibliographystyle{plain} 
\bibliography{refs} 
\end{document}